\documentclass[10pt]{article}
\usepackage{amssymb,amsfonts,amsmath}
\usepackage{enumerate,float, textcomp}
\usepackage{color}
\usepackage{tikz}\usetikzlibrary{tikzmark}\usetikzlibrary{calc}
\usepackage{tcolorbox}
\usepackage{pict2e}
\usepackage{authblk}
\usepackage{xcolor}
\usepackage{hyperref}
\usepackage{dsfont}
\usetikzlibrary{arrows,decorations,backgrounds}
\usetikzlibrary{calc, knots, hobby, arrows.meta}
\pgfdeclarelayer{background layer}
\pgfdeclarelayer{foreground layer}
\pgfsetlayers{background layer,main,foreground layer}
\tikzstyle{arrow} = [thick,-{Stealth[length=5mm]}]

\usepackage[natbib=true, backend = bibtex, style=numeric-comp, sorting=none]{biblatex}
\usepackage{amsthm}

\def\be{\begin{eqnarray}}
\def\ee{\end{eqnarray}}
\def\nn{\nonumber}

\def\tr{{\rm tr}\,}
\def\Tr{{\rm Tr}\,}

\definecolor{red}{rgb}{1,0,0}
\definecolor{orange}{rgb}{1,0.5,0}
\definecolor{violet}{rgb}{0.7,0,1}

\newcommand{\bbl}{{{-\hspace{-1.5mm}\bigcirc\hspace{-1.5mm}-}}}
\newcommand{\thg}{{{\bigcirc \hspace{-12.5pt}-\hspace{-7.5pt}-}}}

\hoffset= - 1.0in      

\begin{document}

\title{\bf Diagrammatic technique for Vogel's universality}

\author[1]{{\bf D. Khudoteplov}\thanks{\href{mailto:khudoteplov.dn@phystech.edu}{khudoteplov.dn@phystech.edu}}}
\author[1,2,3]{{\bf A. Sleptsov}\thanks{\href{mailto:sleptsov@itep.ru}{sleptsov@itep.ru}}}

\vspace{4.5cm}

\affil[1]{Moscow Institute of Physics and Technology, 141700, Dolgoprudny, Russia}
\affil[2]{Institute for Information Transmission Problems, 127051, Moscow, Russia}
\affil[3]{NRC "Kurchatov Institute", 123182, Moscow, Russia\footnote{former Institute for Theoretical and Experimental Physics, 117218, Moscow, Russia}}
\renewcommand\Affilfont{\itshape\small}

\date{}

\maketitle

\vspace{-6.8cm}

\begin{center}
	\hfill MIPT/TH-15/26\\
	\hfill ITEP/TH-17/26\\
	\hfill IITP/TH-15/26
\end{center}

\vspace{4.2cm}

\begin{abstract}
{In his 1999 preprint "Universal Lie Algebra", P. Vogel put forward a hypothesis on the existence of a universal Lie algebra. Although this hypothesis remains open, it is known that many quantities in Lie theory admit universal descriptions. Remarkably, almost all such universal formulas have been obtained through the representation theory of simple Lie (super)algebras, whereas Vogel's original framework was based on a more abstract diagrammatic algebra. Nevertheless, the diagrammatic approach has received little attention over the past two decades, since the last contributions by P.Vogel and J.Kneissler.

In this work, we revive the diagrammatic technique grounded in Vogel's $\Lambda$-algebra and show that it enables truly universal computations. We examine numerous examples and discuss them.
	}
\end{abstract}

\section{Introduction}
\label{sec:introduction}

The study of Lie algebras goes back to the 19th century, and many results have been achieved since then, including classification of semisimple algebras and the numerous applications ranging from knot theory to the particle physics. Ever since the classification of simple Lie algebras was developed, it has been a goal for mathematicians to unveil some universal structure behind them, especially for the exceptional Lie algebras. At the end of the previous century an intriguing discovery was made. Pierre Deligne observed that all of the exceptional Lie algebras as well as the \( \mathfrak{sl}(2)\), \(\mathfrak{sl}(3)\) and \( \mathfrak{so}(8)\) were organised into a 1--parameteric family that was named an exceptional Lie algebra series. At the same time, Pierre Vogel  developed a generalisation that incorporates all the simple Lie algebras into what is called the \textit{Vogel plane}.

Originally, Vogel was researching the finite type (Vassiliev) invariants in knot theory, which led him to a discovery of a universal structure behind all simple Lie algebras and superalgebras. He introduced an algebra, called \( \Lambda \), consisting of trivalent diagrams with three legs modulo the AS and IHX relations, which resemble Lie algebra axioms. Broadly speaking, due to the numerous relations, \( \Lambda \) goes down to three multiplicative generators \( t\), \( \sigma \) and \( \omega \). Each simple Lie algebra provides a numeric function on \( \Lambda \), called \textit{character}. As a side effect, some Lie-algebraic quantities (e.g. Lie algebra dimension) turned out to be expressable by functions in \( t\), \( \sigma \) and \( \omega \), reproducing all the correct answers for a specific Lie algebra \( L \) when mapped by its character \( \chi_{L} \). This concept is known as \textit{Vogel universality}, although more often it is formulated in terms of a different set of parameters \( \alpha \), \( \beta \) and \( \gamma \), which are a Galois extension of \( t\), \( \sigma \) and \( \omega \). 

This observation suggests that there might be some universal Lie algebra, formulated in terms of diagrams, with \( \alpha, \beta, \gamma \) being its parameter space, so that all of the Lie algebras would appear as some specialization of these parameters. In 1999 Vogel published a preprint \cite{vogel1999universal}, dedicated to this particular goal of constructing such an object. However, the paper still remains unfinished, and the 1999 preprint is outdated, since it contains conjectures that have been disproven.

Vogel's finding sparked interest among many researchers, who sought to express more quantities in terms of Vogel parameters. Concrete results include the universal Casimir eigenvalues \cite{mkrtchyan2012casimir, isaev2022split}, quantum dimension \cite{mkrtchyan2017universal}, Chern-Simons partition function \cite{mkrtchyan2012universality, mkrtchyan2013nonperturbative, krefl2013refined}, volume of groups \cite{khudaverdian2017universal}, quantum knot polynomials \cite{mironov2016universal}, Racah matrices (6j-symbols) \cite{mironov2016universalA}. Universal decompositions of adjoint tensor powers have been carried out up to the 5th power \cite{landsberg2006universal, avetisyan2024uniform, isaev2024split}.  It must be noted that all these papers do not utilize the diagrammatic framework, instead, they make use of various methods, specific to each particular problem.

In contrast, some of the formulae turned out to be (at least partially) incompatible with Vogel's universality. For example, at the present moment there is no beautiful universal formula for the rank of Lie algebras. Such examples include the partition function for the refined Chern-Simons theory \cite{krefl2013refined}, Macdonald dimension \cite{bishler2025refined,bishler2025macdonald,bishler2025vogel}. In this cases, only the simple-laced algebras admit universality, which is contrary to the whole Vogel's framework not detecting the difference between even and odd orthogonal algebras.

In this paper we restore to the diagrammatic language in order to systematize the existing knowledge of Vogel universality. We aim to show that the diagrammatic technique is an intuitive framework for universal calculations, capable of producing universal formulae on its own. Also, we demonstrate that the diagrammatic technique is not simply a repackaging of known results, since there exist diagrams, undetected by all semisimple Lie algebra weight systems.

This paper is organised as follows. In Section~\ref{sec:diagrams} we introduce the diagrammatic language, providing definitions of the key terms, such as Jacobi diagrams and Lie algebra weight system. Section~\ref{sec:Vogel} is dedicated to the Vogel framework: we review the \( \Lambda\) algebra, its characters, and the Vogel parameterization. In Section~\ref{sec:kernel} we review the kernel of Lie algebra weight systems, emphasizing the importance of the diagrammatic approach. Section~\ref{sec:universal-formulae} contains a set of universal calculations for different applications. Section~\ref{sec:discussion} contains a brief discussion of an attempt to extend the Vogel's universality from the locus of Lie algebras to the modules on diagrams.

\bigskip

\section{Diagrammatic technique}
\label{sec:diagrams}

\subsection{Motivation}
\label{sec:motivation}
In this section we introduce the Lie algebra diagrammatic technique. Starting from the Feynman diagrams for the Yang-Mills theory, we describe the main concepts of the technique that is used in this paper.

In the Yang-Mills theory the amplitudes are represented by the well-known Feynman diagrams that consist of the 3 types of edges, corresponding to 3 types of particles. The quarks are colored by some representation of the gauge group and the corresponding edge is depicted by a bold line equipped with an orientation representing the flow of charge. Gluons (gauge bosons) are considered in the adjoint representation and represented by thin lines. The ghost field is also in the adjoint representation. 

The collection of the Feynman rules are given in \eqref{eq:YM-boson-fermion}, \eqref{eq:YM-3boson}:

\begin{equation}\label{eq:YM-boson-fermion}
	\mbox{\tikz[baseline = (current bounding box.center)]{
			\coordinate (a) at (0,0);
			\node (b) at (0,1) {\( a, \mu \) };
			\draw (a)--(b);
			\draw[ultra thick] (a)--+(210:1);
			\draw[-{Stealth[length=3mm]}, ultra thick] (a) --+(330:1);
			\node at ($(a)+(2,0)$) {\large \( =i g \gamma^{\mu} \color{red}T^a \) };
	}}
\end{equation}
\begin{equation}\label{eq:YM-3boson}
	\mbox{\tikz[baseline = (current bounding box.center)]{
			\coordinate (o) at (0,0);
			\node (a) at (90:1.2) {\( a, \mu \) };
			\node (b) at (210:1.5) {\( b, \nu \) };
			\node (c) at (330:1.5) {\( c, \rho \) };
			\draw (a) -- (o) --(c);
			\draw (b)--(o);
			\draw[-{Stealth[length=1.5mm]}, thick] (0.15, 0.7)--+(270:0.5);
			\draw[-{Stealth[length=1.5mm]},  thick] (195:0.7)--+(30:0.5);
			\draw[-{Stealth[length=1.5mm]},  thick] (315:0.7)--+(150:0.5);
			\node at (0.33, 0.5) {\( k \) };
			\node at (165:0.55) {\( p \) };
			\node at (285:0.6) {\( q \) };
			\node at (5.,0) {\( =g {\color{red}f^{abc}} [g^{\mu \nu} (k-p)^\rho+g^{\nu \rho} (p-q)^\mu+g^{\rho \mu}(q-k)^\nu] \) };
	}}
\end{equation}

The amplitudes can be factorised in two parts: the kinematic one and the color one. In this paper we focus on the latter, considering only the group-theoretical part of the Feynman rules. An accurate definition of our framework is given in the rest of this section.

\subsection{Collection of rules}
\label{sec:rules}

We work with the trivalent diagrams modulo the so-called STU, AS and IHX relations, depicted in \eqref{eq:STU-relation} and \eqref{eq:AS-IHX}. Thick lines are oriented and the inner vertices have their half-edges cyclic ordered (in this paper we assume the counterclockwise orientation). In this paper we call them Jacobi diagrams. We differentiate between several types of Jacobi diagrams, using the following terms:

\begin{itemize}
	\item Closed Jacobi diagrams \( \mathcal{A}\)~--- these are diagrams that have no loose legs and necessarily involve a Wilson loop
	\item Open Jacobi diagrams \( \mathcal{B} \)~--- these diagrams have no Wilson lines and symmetrized over its legs
	\item Primitive Jacobi diagrams \( \mathcal{P}\)~--- these are closed Jacobi diagrams with connected inner graph or open Jacobi diagrams with one connected component
	\item Chord diagrams --- closed Jacobi diagrams that have no inner vertices
\end{itemize}

STU  relation mean that Lie algebra representation represent Lie bracket as a commutator, while the AS and IHX identities translate to the antisymmetry of Lie bracket and the Jacobi identity.

\begin{equation}
	\label{eq:STU-relation}
	\mbox{\tikz[scale=0.5, baseline = (current bounding box.center)]{
			\coordinate (s) at (0,0);
			\draw[-{Stealth[length=3mm]}, ultra thick] ($(s)+ (45:-2)$) arc(225:315:2);
			\coordinate (t) at (4,0);
			\draw[-{Stealth[length=3mm]}, ultra thick] ($(t)+ (45:-2)$) arc(225:315:2);
			\coordinate (u) at (8,0);
			\draw[-{Stealth[length=3mm]}, ultra thick] ($(u)+ (45:-2)$) arc(225:315:2);
			\coordinate (s1) at ($(s)+(-1, 0)$);
			\coordinate (t1) at ($(t)+(-1, 0)$);
			\coordinate (u1) at ($(u)+(-1, 0)$);
			\coordinate (s2) at ($(s)+(1, 0)$);
			\coordinate (t2) at ($(t)+(1, 0)$);
			\coordinate (u2) at ($(u)+(1, 0)$);
			\draw ($(s)+(0, -2)$) -- ($(s)+(0, -1)$);
			\draw (s1) -- ($(s)+(0, -1)$);
			\draw ($(s)+(0, -1)$) -- (s2);
			\draw[fill] ($(s)+(0, -1)$) circle(0.1);
			\draw (t1) -- ($(t)+(70:-2)$);
			\draw ($(t)+(110:-2)$) -- (t2);
			\draw (u1) -- ($(u)+(110:-2)$);
			\draw ($(u)+(70:-2)$) -- (u2);
			\node[scale=1.5] at ($0.5*(s) + 0.5*(t)-(0,1)$) {\( - \) };
			\node[scale=1.5] at ($0.5*(u) + 0.5*(t)-(0, 1)$) {\( + \) };
			\node[scale=1.5] at ($0.5*(s) - 0.9*(t)-(0, 1)$) {{\footnotesize \textbf{STU}}  \  $=$};
			\node[scale=1.5] at ($ (u) + (2.5,-1) $) {$ = 0 $};
	}}
\end{equation}

\begin{equation}\label{eq:AS-IHX}
	\mbox{\begin{picture}(850,20)(-60,30)
			\put(0,0){\line(0,1){20}}
			\qbezier(0,20)(10,25)(-10,45)
			\qbezier(0,20)(-10,25)(10,45)
			\put(20,17){\Large$=   \ -$}
			\put(60,0){\line(0,1){20}}
			\put(60,20){\line(-1,2){12.5}}
			\put(60,20){\line(1,2){12.5}}
			\put(20,-12){\mbox{\fontsize{12}{12}$\textbf{AS}$}}
			\put(170,10){\line(0,1){25}}
			\put(170,35){\line(-2,1){20}}
			\put(170,35){\line(2,1){20}}
			\put(170,10){\line(-2,-1){20}}
			\put(170,10){\line(2,-1){20}}
			\put(210,17){\Large $=$}
			\put(238,44){\line(1,-2){12}}
			\put(240,0){\line(1,2){10}}
			\put(250,20){\line(1,0){25}}
			\put(287,44){\line(-1,-2){12}}
			\put(285,0){\line(-1,2){10}}
			\put(305,17){\Large $-$}
			\put(330,0){\line(1,2){10}}
			\put(330,43.5){\line(3,-2){35}}
			\put(340,20){\line(1,0){25}}
			\put(375,0){\line(-1,2){10}}
			\put(340,20){\line(3,2){35}}
			\put(250,-12){\mbox{\fontsize{12}{12}$\textbf{IHX}$}}
	\end{picture}}
	\vspace{15mm}
\end{equation}

Trivalent vertices correspond to the structure constants, and the thick oriented lines (they are typically called Wilson lines/loops in this paper) correspond to some representation or to the universal enveloping algebra. Indices are contracted using a fixed nondegenerate invariant metric.

No lie algebra is specialized, hence, there are no rules apart from the aforementioned.
Note that there is no ordering of the half-edges for the vertices on the Wilson line. Instead, the Wilson lines are oriented, much like in the ordinary QFT. 

One can specify the Wilson lines to be colored by the adjoint representation. This mapping is denoted by $ \pi_\text{adj} $ in this paper:

\begin{equation}\label{eq:Wilson-to-adjoint}
	\mbox{\tikz[use Hobby shortcut, baseline = (current bounding box.center)]{
			\coordinate (d) at (-4,0);
			\coordinate (a) at (0,0);
			\coordinate (b) at (4,0);
			\coordinate (c) at (8,0);
			\node at ($0.5*(a)+0.5*(d)+(0, 0.4)$) {\huge$=$};
			\node at ($0.5*(a)+0.5*(b)+(0, 0.4)$) {\huge$=$};
			\node at ($0.5*(b)+0.5*(c)+(0, 0.4)$) {\huge$=$};
			\node at ($(d)+(-1.2,0.4)$) {\Huge (};
			\node at ($(d)+(-1.7,0.3)$) {\large $ \pi_\text{adj}$};
			\node at ($(d)+(1.2,0.4)$) {\Huge )};
			\draw[-{Stealth[length=3mm]}, ultra thick] ($(d)+(-1,0)$) -- ($(d)+(1,0)$);
			\draw (d)-- +(0,1);
			\draw[-{Stealth[length=3mm]}, ultra thick] ($(a)+(-1,0)$) -- ($(a)+(1,0)$);
			\draw (a)-- +(0,1);
			\node at ($(a)+(-0.8, -0.2)$) {\small adj}; 
			\draw[-{Stealth[length=3mm]}, ultra thick] ($(b)+(-1,0)$) -- ($(b)+(1,0)$);
			\draw (b) .. +(0,-0.05)..+(0.05,-0.1) .. +(0.15, 0) .. +(0.1, 0.3) .. + (0.04, 0.7) .. +(0, 1);
			\node at ($(b)+(-0.8, -0.2)$) {\small adj}; 
			\draw ($(c)+(-1,0)$) -- ($(c)+(1,0)$);
			\draw (c)-- +(0,1);
	}}
\end{equation}
The second equality sign in \eqref{eq:Wilson-to-adjoint} is to illustrate that there is no ordering of the half-edges on the Wilson line.

\subsection{Weight system / Lie algebra specialization}
\label{sec:weght-system}

Given a diagram modulo the AS, IHX and STU relations, one can specialize a particular Lie algebra. The structure constant $f^{abc}$ is assigned to each inner vertex and the metric tensor $g_{ab}$ is assigned to each edge. Outer vertices (on the Wilson loops) correspond to the Lie algebra generators taken in some representation. Then the duplicate indices are contracted, giving rise to an element of $L^{\otimes n} \otimes_i \text{End}(R_i) $, where $n$ is the number of legs of the diagram that correspond to free uncontracted tensor indices and $i$ labels the Wilson lines, with $ R_i $ being the representation for  the $i$-th Wilson line. 

This mapping is called a Lie algebra weight system and denoted by $\Phi_L$ (if there are Wilson lines, representations are also specified), where $L$ is a pair of Lie algebra and a non-degenerate invariant metric (typically the Killing form).

\begin{equation}
	\large\Phi_L\Big( 
	\mbox{\tikz[ scale = 2, baseline=8pt]{
			\coordinate (a) at (0,0);
			\draw [ultra thick, -{Stealth[length=2mm]}] (a) -- +(1,0);
			\draw (a)+(0.5, 0) arc(0:180:0.2);
			\draw (a)++(0.7, 0) -- ++(0, 0.2)-- + (-0.15, 0.2);
			\draw (a)++(0.7, 0) -- ++(0, 0.2)-- + (0.15, 0.2);
	}} \
	\Big) = T_i T^i f^{abc} T_c
\end{equation}

Diagrams have proven to be a useful tool for calculations in Lie algebras and their representation theory. We recommend \cite{cvitanovic2008group} as a comprehensive text on this topic. In the next Section we focus on the Vogel framework, which is specific to the adjoint sector of representation theory.

\section{Vogel framework}
\label{sec:Vogel}

\subsection{Structure and basis}
\label{sec:Lambda-structure}
$\Lambda$ is a vector space generated by 3-legged diagrams of valence 3 modulo the AS and IHX relations. There is also a condition for a diagram in $\Lambda$ to be antisymmetric with respect to permutation of its legs, however, it is unknown whether there are any 3-legged diagrams other than antisymmetric.

There exists a way to multiply any connected trivalent diagram modulo AS+IHX by an element of $\Lambda$. To do so, one should  insert a 3-legged $\Lambda$-diagram into a trivalent vertex of the other diagram. Hence, $\Lambda$ is an algebra with respect to that vertex multiplication and some other spaces of diagrams (like primitive Jacobi diagrams $\mathcal{P}$) acquire the structure of $\Lambda$-module:

\begin{equation}\label{eq:lambda-product}
	\mbox{\tikz[scale=0.7, baseline = (current bounding box.center)]{
			\node[inner sep=0.5mm] (a) at (0.2,-0.1) {$\hat{u}$};
			\coordinate (b) at (3,0);
			\node[inner sep=0.5mm] (c) at (6.4,0) {$\hat{u}$};
			\draw[fill] (a)+(1.2,0.05) circle(0.05);
			\node at ($0.5*(c)+0.5*(b)$) { \(=\)};
			\draw (a) -- +(-30:0.8);	
			\draw (a) -- +(90:0.8);	
			\draw (a) -- +(210:0.8);	
			\draw[ultra thick] (b) circle(1);
			\draw (b) -- + (0:1);
			\draw (b) -- + (120:1);
			\draw (b) -- + (240:1);
			\draw[fill] (b) circle(0.04);
			\draw[ultra thick] (c) circle(1);
			\draw (c) -- + (0:1);
			\draw (c) -- + (120:1);
			\draw (c) -- + (240:1);
	}}
\end{equation}

$\Lambda$-algebra was extensively studied by P. Vogel \cite{vogel2011algebraic} and Jan Kneissler \cite{kneissler2001spaces}. However, there are still several gaps that remain unknown. The current state of the art goes down to the following:
\begin{enumerate}
	\item $\Lambda$ is associative and commutative algebra.
	\item \( \Lambda \) is graded by half the number of trivalent vertices minus one half.
	\item All the known elements in $\Lambda$ are generated by $\hat{t}= \mbox{\tikz[scale=0.60, baseline = -1.5pt]{
		\draw (0,0) circle(0.3); 
\draw (90:0.3)--(90:0.7); 
\draw (210:0.3)--(210:0.7); 
\draw (330:0.3)--(330:0.7); 
}}$ 
and $\hat{x}_n = \mbox{\tikz[scale=0.40, baseline = -1pt]{
\draw (0,0.6) --(0,1.25);
\draw (-0.45,0.6) arc(90:270:0.3);
\draw (0.45,0.6) arc(90:-90:0.3);
\draw (-0.45,0.6) -- (0.45,0.6);
\draw (-0.5,0) -- (0.5,0);
\draw (-0.45,0) -- +(0,-0.5);
\draw (-0.2,0) -- +(0,-0.5);
\draw (0.45,0) -- +(0,-0.5);
\draw (-1.2,-0.5) -- (1.2, -0.5);
\draw[fill] (0, -0.25) circle(0.01);
\draw[fill] (0.125, -0.25) circle(0.01);
\draw[fill] (0.25, -0.25) circle(0.01);
}}$ diagrams, but there are relations between them.
\item The relations become simple in another basis of the $\hat{\omega}_p$ diagrams: $\hat{\omega}_p \hat{\omega}_q = \hat{\omega}_0 \hat{\omega}_{p+q}$. This allows to introduce an auxiliary variable \( \sigma \) so that \( \sigma^p \hat\omega_0 = \hat\omega_p \). 
	\item $ \Lambda $ has zero divisors.
\end{enumerate}
\noindent \textit{Remark.} We use variables with hats for elements of $\Lambda$-algebra. Note that $\sigma$ is not the element of $\Lambda$ in contrast to $\sigma \hat{\omega}$, where $\omega\equiv \omega_0$.

\subsection{Characters}
\label{sec:char}

The correspondence between Lie algebra structure and the AS+IHX relations was described in Section~\ref{sec:weght-system}. For an element of $\Lambda$ structure constant $f^{abc}$ is assigned to each vertex and metric tensor $g_{ab}$ (non-degenerate Killing form) is assigned to each edge:
\begin{align}
\Phi_{L}: \Lambda \rightarrow {\textstyle\bigwedge}^3 L, \
		\put(23,3){\line(0,1){15}}
		\put(23,3){\line(3,-2){14.5}}
		\put(23,3){\line(-3,-2){14.5}}
		\put(23,3){\circle*{2}}
		\hspace{12mm} \mapsto \ f^{abc}
		\nonumber
	\end{align}

After contraction of the duplicate indices we end up with a tensor with 3 free indices. But for simple Lie algebra there is only one invariant antisymmetric tensor of rank 3, namely the structure constant:
\begin{align}
	\label{eq:char}
	\Phi_{L}(\hat{x}) = \chi_{L}(\hat{x}) \, f^{abc}, \qquad \forall\, \hat{x} \in \Lambda.
\end{align}

The coefficient of proportionality between the resulting tensor and the structure constant is called a character of the Lie algebra taken on that element of $\Lambda$. The character on $\Lambda$ coming from a Lie algebra $L$ defined as above is denoted by $\chi_{L}: \Lambda \rightarrow \mathbb{Q}$.

It is easy to check that character is a ring homomorphism. 
Let $\chi_{L}(\hat{t})=t_L, \ \chi_{L}(\hat{\omega})=\omega_L$ and $\chi_{L}(\sigma\hat{\omega})=\sigma_L\omega_L$. From now we basically omit the index $L$ in variables $t, \sigma, \omega$ to simplify notations. Hence, for the elements of $\Lambda$ expressible in the form of a polynomial in $t$, $\sigma$ and $\omega$ (currently, all known elements of $\Lambda$ are of this type) one only needs to know the values that $\chi_{L}$ takes on $t$, $\sigma$ and $\omega$. Those depend on ${L}$ and on the choice of metric in $L$. For any renormalization of the metric $t$, $\sigma$ and $\omega$ also rescale as variables of degree 1, 2 and 3 correspondingly. Hence, the values $\chi_{L}(\hat{\omega})/\chi_{L}(\hat{t}^3)$ and $\chi_{L}(\sigma\hat{\omega})/\chi_{L}(\hat{t}^2\hat{\omega})$  are invariant under rescaling and depend exclusively on the Lie algebra structure.

Instead of variables $t, \sigma, \omega$ there is another rather convenient parametrization via $\alpha, \beta, \gamma$:
\begin{equation}
	\label{eq:alphaparam}
	\alpha + \beta +\gamma = t, \quad \alpha \beta + \beta \gamma + \alpha \gamma = \sigma-2t^2, \quad \alpha \beta \gamma = \omega-t\sigma\,.
\end{equation}

In terms of $\alpha,\ \beta,\ \gamma$ the table of characters is shown in Table \ref{tab:char} in the same normalization as in \cite{mkrtchyan2012casimir}:
\renewcommand{\arraystretch}{1.5}
\begin{table}[h]
	\centering
	\begin{tabular}{|c|c|c|c|}
		\hline
		Lie algebra $L$& $\alpha$ & $\beta$ & $\gamma$ \\
		\hline
		$\mathfrak{sl}_{N}$ & $-2$ & $2$ &$N$ \\
		\hline
		$\mathfrak{so}_{N}$ & $-2$ & $4$ & $N-4$\\
		\hline
		$\mathfrak{sp}_{2N}$ & $-2$ & $1$ &$N+2$ \\
		\hline
		$\mathfrak{g}_2$ & $-2$ & $10/3$ & $8/3$ \\
		\hline
		$\mathfrak{f}_4$ & $-2$ & $5$ &$6$ \\
		\hline
		$\mathfrak{e}_6$ & $-2$ & $6$ & $8$ \\
		\hline
		$\mathfrak{e}_7$ & $-2$ & $8$ & $12$\\
		\hline
		$\mathfrak{e}_8$ & $-2$ & $12$ & $20$ \\
		\hline
	\end{tabular}
	\caption{Vogel parameters.}	
	\label{tab:char}
\end{table}

\bigskip

	\subsection{Zero divisor} 
	\label{sec:zero-divisor}
	We have already mentioned that $\Lambda$-algebra has a divisor of zero $\hat{t} \cdot \hat{P}_{15} = 0$. This is proven in \cite[Theorem 8.4.]{vogel2011algebraic} and now we briefly describe the construction. Let ${U}$ be the following 6-legged diagram:
	
	\mbox
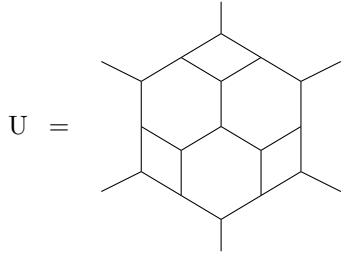
{\begin{picture}(300,120)(-150,-5)
	\put(0,50){${\text{U}} \ \, = $}
	\put(50,53){\line(0,1){17}}
	\put(80,53){\line(0,1){17}}
	\put(50,70){\line(5,3){15}}
	\put(65,79){\line(5,-3){15}}
	\put(50,53){\line(5,-3){15}}
	\put(110,53){\line(0,1){17}}
	\put(80,70){\line(5,3){15}}
	\put(95,79){\line(5,-3){15}}
	\put(95,44){\line(5,3){15}}
	\put(65,27){\line(0,1){17}}
	\put(95,27){\line(0,1){17}}
	\put(65,44){\line(5,3){15}}
	\put(80,53){\line(5,-3){15}}
	\put(80,18){\line(5,3){15}}
	\put(65,27){\line(5,-3){15}}
	\put(50,36){\line(0,1){17}}
	\put(50,36){\line(5,-3){15}}
	\put(110,36){\line(0,1){17}}
	\put(95,27){\line(5,3){15}}
	\put(80,88){\line(5,-3){15}}
	\put(65,79){\line(5,3){15}}
	\put(80,88){\line(0,1){12}}
	\put(80,6){\line(0,1){12}}
	\put(50,70){\line(-2,1){15}}
	\put(50,36){\line(-2,-1){15}}
	\put(110,36){\line(2,-1){15}}
	\put(110,70){\line(2,1){15}}
	\end{picture}}
	
	Then we  \text{define}  $\hat{P}_{15} \neq 0 \in \Lambda$ by  removing a neighborhood of a trivalent  vertex from $V$:
	
	\mbox{\begin{picture}(850,30)(-375,-25)
	\put(-315,-22){V= {\Large  \text{$\frac{1}{6!}\sum \limits_{\sigma \in S_6}$}} $ \text{sign}(\sigma) \ \cdot$}
	\put(-199,-20){\circle{30}}
	\put(-203,-22){${\text{U}}$}
	\put(-197,-35){\line(1,0){32}}
	\put(-187,-29){\line(1,0){22}}
	\put(-184.5,-23){\line(1,0){19.5}}
	\put(-184.5,-17){\line(1,0){19.5}}
	\put(-187,-11){\line(1,0){22}}
	\put(-197,-5){\line(1,0){32}}
	\put(-165,-40){\line(0,1){40}}
	\put(-165,-40){\line(1,0){20}}
	\put(-165,-0){\line(1,0){20}}
	\put(-145,-40){\line(0,1){40}}
	\put(-157,-22){$\sigma$}
	\put(-145,-35){\line(1,0){32}}
	\put(-145,-29){\line(1,0){22}}
	\put(-145,-23){\line(1,0){19.5}}
	\put(-145,-17){\line(1,0){19.5}}
	\put(-145,-11){\line(1,0){22}}
	\put(-145,-5){\line(1,0){32}}
	\put(-111,-20){\circle{30}}
	\put(-114,-22){${\text{U}}$}
	\put(-90,-22){$=$}
	\put(-60.7,-20){\circle{25}}
	\put(-68,-22){$\hat{P}_{15}$}
	\put(-48,-20){\line(1,0){14.5}}
	\put(-53.3,-10){\line(2,-1){20}}
	\put(-53.3,-30){\line(2,1){20}}
	\put(-34,-20){\circle*{2}}
	\qbezier[8](-34,-15)(-29,-15)(-29,-20)
	\qbezier[8](-29,-20)(-29,-25)(-34,-25)
	\qbezier[8](-34,-25)(-39,-25)(-39,-20)
	\qbezier[8](-39,-20)(-39,-15)(-34,-15)
	\put(-20,-15){\footnotesize remove}
	\put(-21,-15){\vector(-2,-1){10}}
	\end{picture}}
	
	\vspace{0mm}
	\begin{align}
&\omega P_{\mathfrak{sl}}P_{\mathfrak{osp}}P_{\mathfrak{exc}} = -2^{10} \hat{P}_{15} \nonumber \\
&\hat{t} \cdot \hat{P}_{15} = 0 \nonumber, \quad \hat{t} = 
\mbox{\put(10,0){\line(1,0){16}}
	\put(10,0){\line(1,2){8}}
	\put(18,16){\line(1,-2){8}}
	\put(18,16){\line(0,1){8}}
	\put(10,0){\line(-2,-1){8}}
	\put(26,0){\line(2,-1){8}}}
	\end{align}

Symmetric polynomials of $\alpha,\beta,\gamma$ which vanish for particular sets of algebras in  (\ref{tab:char}) are:
\be
\label{Palgebra0}
P_{\mathfrak{sl}} &=& (\alpha+\beta)\,(\beta+\gamma)\,(\alpha+\gamma), \hspace{0mm} \nonumber \\
P_{\mathfrak{osp}} &=& (\alpha+2\beta)(2\alpha+\beta)\,(\beta+2\gamma)(2\beta+\gamma)\,(\alpha+2\gamma)(2\alpha+\gamma),   \\
P_{\mathfrak{exc}} &=& (\alpha-2\beta-2\gamma)\,(\beta-2\alpha-2\gamma)\,(\gamma-2\alpha-2\beta). \nonumber
\ee
See more details about zero divisor in $\Lambda$-algebra in \cite{vogel2011algebraic} and in \cite{khudoteplov2025yangbaxter}.

	\subsection{0-,1- and 2-legged trivalent diagrams}
	\label{sec:0-1-2-legs}
	
	\paragraph{0-legged diagrams.}
	
	Trivalent diagrams (without external legs) modulo the AS and IHX relations are also called 3-graphs in literature \cite{chmutov1998algebra}. The space of 3-graphs is graded by half the number of vertices and is denoted by $\Gamma = \bigoplus_{n=0}^{\infty} \Gamma_n$. The space of 3-graphs $\Gamma$ admits two different products, the vertex product and the edge product. As an algebra with the vertex product, the 3-graphs with at least one vertex $\Gamma_{n \geq 1}$ are isomorphic to Vogel's $ \Lambda $ (see Chapter 7 in \cite{chmutov2012introduction}). The mapping between $\Gamma_{n \geq 1}$ and $ \Lambda $ is the following:

	\begin{equation}
\label{eq:3graphs-Lambda}
\mbox{\tikz[use Hobby shortcut, scale=0.6, baseline = (current bounding box.center)]{
		\coordinate (a) at (0, 0);
		\coordinate (b) at (5, 0);
		\draw[fill, opacity=0.2] (a) ellipse(1.5 and 1);
		\draw[fill, opacity=0.2] (b) ellipse(1.5 and 1);
		\node at ($ 0.5*(a)+0.5*(b)+(0,-0.2) $) {\Large \(  \mapsto  \) };		
		\node at ($(a)+(-3,0)$) {\Large $\Xi :$};
		\draw (a) ++(-1,-0.75) .. ++(-0.1,-0.2) .. ++(0, -0.1) .. ($ (a)+(0, -2) $);
		\draw (a) ++(1,-0.75) .. ++(0.1,-0.2) .. ++(0, -0.2) .. ($ (a)+(0, -2) $);
		\draw (a) ++ (0, -1) .. ++ (0, -1);
		\draw [fill] ($ (a)+(0, -2) $) circle(0.03);
		\draw (b) ++(-1,-0.75) .. ++(-0.05,-0.2) .. ++(0, -0.1) .. ($ (b)+(-1, -1.2) $);
		\draw (b) ++(1,-0.75) .. ++(0.05,-0.2) .. ++(0, -0.1) .. ($ (b)+(+1, -1.2) $);
		\draw (b) ++ (0, -1) .. ++ (0, -0.2);
		\draw ($ (b)+(-1.2, -1.7) $) rectangle($ (b)+(+1.2, -1.2) $);
		\node at ($ (b)+(0, -1.45) $) {\small antisym};
		\draw ($ (b)+(0, -1.7) $)--  ++ (0, -0.3);
		\draw ($ (b)+(-.9, -1.7) $)--  ++ (0, -0.3);
		\draw ($ (b)+(0.9, -1.7) $)--  ++ (0, -0.3);
}}
\end{equation}

\noindent The gray zone is a subgraph modulo the AS and IHX relations that is the same on the both sides of equation.
For example:

\begin{align}\mbox{
	\put(0,3){\circle{20}}
	\put(-10,3){\circle*{2}}
	\put(10,3){\circle*{2}}
	\put(45,3){\circle*{2}}
	\put(-10,3){\line(1,0){20}}
	\put(20,0){$\mapsto$}
	\put(45,3){\line(0,1){10}}
	\put(45,3){\line(3,-2){9}}
	\put(45,3){\line(-3,-2){9}}
} \quad \hspace{17mm}, \quad\quad \mbox{\tikz[scale=0.3, baseline, anchor=base]{
		\coordinate (a) at (0,0.5);
		\coordinate (b) at (4,0.5);
		\coordinate (c) at (8, 0.7);
		\draw (a) circle(1);
		\draw (b) circle(1);
		\draw ($(a)+(10:1)$)--($(b)+(170:1)$);
		\draw ($(a)+(30:1)$)--($(b)+(150:1)$);
		\draw ($(a)+(-10:1)$)--($(b)+(190:1)$);
		\draw ($(a)+(-30:1)$)--($(b)+(210:1)$);
		\draw[fill] ($(a)+(10:1)$) circle(0.05);
		\draw[fill] ($(a)+(30:1)$) circle(0.05);
		\draw[fill] ($(a)+(-10:1)$) circle(0.05);
		\draw[fill] ($(a)+(-30:1)$) circle(0.05);
		\draw[fill] ($(b)+(170:1)$) circle(0.05);
		\draw[fill] ($(b)+(150:1)$) circle(0.05);
		\draw[fill] ($(b)+(190:1)$) circle(0.05);
		\draw[fill] ($(b)+(210:1)$) circle(0.05);
		\draw (c) circle(0.7);
		\draw (c)+(-0.9, -1.5) -- +(0.9, -1.5);
		\draw ($(c)+(60 :-0.7)$) -- ($(c)+(-0.5, -1.5)$);
		\draw ($(c)+(90 :-0.7)$) -- ($(c)+(0, -1.5)$);
		\draw ($(c)+(120 :-0.7)$) -- ($(c)+(0.5, -1.5)$);
		\draw (c)++(0, 0.7) -- + (0, 0.4);
		\node at ($(b)+(2,-0.3)$) {$\mapsto$};
}}
\end{align}

\noindent The only difference between $\Gamma$ and $\Lambda$ is in the one dimensional subspace $\Gamma_0$ generated by the circle $\bigcirc$ without vertices.

\paragraph{1-legged diagrams.}

The case of 1-legged diagrams is the most trivial one. Using the AS and IHX relations, any diagram with $ 1 $ leg can be transformed into a diagram with an isolated loop that vanishes by the AS relation.

\begin{equation}\label{eq:1-legged-zero}
\centering
\mbox{\tikz[use Hobby shortcut, scale=0.6, baseline = (current bounding box.center)]{
		\coordinate (a) at (0,0);
		\coordinate (b) at (5, 0);
		\draw[fill, opacity=0.2] (a) ellipse(1.5 and 1);
		\draw[fill, opacity=0.2] (b) ellipse(1.5 and 1);
		\node at ($ 0.5*(a)+0.5*(b)+(0, -0.1) $) {\Large \(  = \) };	
		\draw (a) ++(-0.6,-0.9) .. ++(0,-0.5) .. ++(0, -0.6);
		\draw (a) ++(0.6,-0.9) .. ++(-0.1,-0.2) .. ++(-0.2, -0.2) .. ++(-0.9, -0.2);
		\draw (b) ++(-0.6,-0.9) .. ++(0,-0.5) .. ++(0, -0.6);
		\draw (b) ++(0.6,-0.9) -- ($ (b)+(0.6, -1.2) $);
		\draw (b)+(0.6, -1.5) circle (0.3);
		\node at ($ (b)+(2.8, 0) $) {\Large \(  = 0  \) };
}}
\end{equation}

\paragraph{2-legged diagrams.}

The space of diagrams with $ 2 $ legs is isomorphic to the aforementioned space of diagrams with $ 0 $ legs. This isomorphism is given by cutting an edge of a 0-legged diagram and, reversely, by gluing the legs of a 2-legged diagram:
\begin{figure}[H]
\centering
\mbox{\tikz[use Hobby shortcut, scale=0.6, baseline = (current bounding box.center)]{
		\coordinate (a) at (0, 0);
		\coordinate (b) at (5, 0);
		\draw[fill, opacity=0.2] (a) ellipse(1.5 and 1);
		\draw[fill, opacity=0.2] (b) ellipse(1.5 and 1);
		\node at ($ 0.5*(a)+0.5*(b) $) {\Large \(  \leftrightarrow  \) };		
		\draw (a) +(-1,-0.75) .. + (-1, -0.8) .. (0, -1.5) .. (1, -0.8) .. (1, -0.75);
		\draw (b) ++(-1,-0.75) .. ++(0,-0.2) .. ++(0, -0.3);
		\draw (b) ++(1,-0.75) .. ++(0,-0.2) .. ++(0, -0.3);
}}
\end{figure}

Since the 0-legged diagrams with at least one vertex are isomorphic to $ \Lambda $, the 2-legged diagrams with at least one vertex are generated as a $ \Lambda $-module by the "bubble" diagram $ \bbl $. Thus it becomes apparent that all 2-legged diagrams are symmetric with respect to the permutation of its legs. 

\bigskip

To sum up, the spaces of diagrams with low number of legs are quite simple. All the diagrams with 1 leg vanish. And the spaces of diagrams with 0 and 2 legs are isomorphic to each other and almost isomorphic to the $\Lambda$, the only difference being the diagrams with no vertices.

Add some pictures as an illustration like
\begin{align}
\put(0,3){\circle{10}}
\put(-5,3){\circle*{2}}
\put(5,3){\circle*{2}}
\put(-15,3){\line(1,0){10}}
\put(5,3){\line(1,0){10}}
\put(20,0){$\leftrightarrow$}
\put(45,3){\circle{20}}
\put(35,3){\circle*{2}}
\put(55,3){\circle*{2}}
\put(35,3){\line(1,0){20}}
\put(65,0){$\mapsto$}
\put(85,3){\line(0,1){10}}
\put(85,3){\line(3,-2){9}}
\put(85,3){\line(-3,-2){9}}
\put(85,3){\circle*{2}}
\put(110,0){$...$}
\end{align}

\bigskip

The aforementioned "bubble" diagram $\bbl$ has a very important Lie-algebraic meaning. Actually, it corresponds to the Killing form, defined as $ K(x, y) = \tr (\text{ad}_x \circ \text{ad}_y)$:

\begin{figure}[h]
\centering
\mbox{\tikz[use Hobby shortcut, scale=0.8, baseline = (current bounding box.center)]{
		\coordinate (a) at (0,0);
		\coordinate (b) at (0, -0.9);
		\draw ($ (a)+(0:1) $) arc (0:180:1);
		\draw ($ (b)+(180:1) $) arc (180:360:1);
		\draw ($ (a) + (-1, 0) $) -- ($ (b) + (-1, 0) $);
		\draw ($ (a) + (1, 0) $) -- ($ (b) + (1, 0) $);
		\node[] (y) at ($ (a)+ (-1.55, 1.6) $) {\small \(  y  \) };
		\node[] (x) at ($ (a)+ (-2.2, 1.6) $) {\small \(  x  \) };
		\draw (a) +(170:1) .. +(-1.2, 0.4) .. + (-1.48, 1) ..($ (y) + (0, -0.15)$)  ; 
		\draw (a) +(-1, -0.8) .. +(-2, 0.4) .. + (-2.15, 1) .. ($ (x) + (0, -0.15)$) ; 
		\draw[dotted] (a)++(-2.5, 1.3) -- ++(4, 0) node[right] {\small\(  x\otimes y  \) }; 
		\draw[dotted] (a)++(-2.5, 0.7) -- ++(4, 0) node[right] {\small\(  x\otimes y \otimes e_i \otimes e^i  \) }; 
		\draw[dotted] (a)++(-2.5, -0.2) -- ++(4, 0) node[right] {\small\(  x\otimes \left[ y, e_i \right] \otimes e^i  \) }; 
		\draw[dotted] (a)++(-2.5, -1.3) -- ++(4, 0) node[right] {\small\(  \left[ x, \left[ y, e_i \right] \right] \otimes e^i  \) }; 
		\draw[dotted] (a)++(-2.5, -2.3) -- ++(4, 0) node[right] {\small\(  K(x, y)  \) }; 
} }
\end{figure}

For a simple Lie algebra there exists only one invariant metric (up to proportionality). Hence, the 'bubble' is proportional to an edge if we apply the mapping $ \Phi_{L} $. The coefficient of proportionality is $ 2 t_{L} $: $ \Phi_{L} (\bbl) = 2 t_{L} \Phi_{L} (-\hspace{-0.6mm}-) $, due to the following equality:

\begin{equation}
\mbox{\tikz[scale=1, baseline=1mm]{
		\coordinate (a) at (0,0);
		\coordinate (b) at (2,0);
		\coordinate (c) at (4,0);
		\coordinate (d) at (6.4,0);
		\node at (1,0.4) {$=$};
		\node at (3,0.4) {$-$};
		\node at (5.25,0.4) {$= \; 2$};
		\draw (a)+(-0.60,0) -- +(0.60,0);
		\draw (b)+(-0.60,0) -- +(0.60,0);
		\draw (c)+(-0.60,0) -- +(0.60,0);
		\draw (d)+(-0.60,0) -- +(0.60,0);
		\draw (a)+(0,0.7) -- +(0,1.0);
		\draw (b)+(0,0.7) -- +(0,1.0);
		\draw (c)+(0,0.7) -- +(0,1.0);
		\draw (d)+(0,0.7) -- +(0,1.0);
		\draw (a)+(0,0.5) circle (0.2);
		\draw (a) -- +(0,0.3);
		\draw (b)+(0,0.7) -- +(-0.3,0);
		\draw (b)+(0,0.7) -- +(+0.3,0);
		\draw (d)+(0,0.7) -- +(-0.3,0);
		\draw (d)+(0,0.7) -- +(+0.3,0);
		\draw ($(c)+(0,0.7)$) to [out=-145,in=135] ($(c)+(0.3, 0)$);
		\draw ($(c)+(0,0.7)$) to [out=-35,in=45] ($(c)+(-0.3, 0)$);
		\draw[fill] (a)+(0,0.7) circle (0.03);
		\draw[fill] (a)+(0,0.3) circle (0.03);
		\draw[fill] (a) circle (0.03);
		\draw[fill] (b)+(0,0.7) circle (0.03);
		\draw[fill] (c)+(0,0.7) circle (0.03);
		\draw[fill] (d)+(0,0.7) circle (0.03);
		\draw[fill] (b)+(0.3,0) circle (0.03);
		\draw[fill] (b)+(-0.3,0) circle (0.03);
		\draw[fill] (c)+(0.3,0) circle (0.03);
		\draw[fill] (c)+(-0.3,0) circle (0.03);
		\draw[fill] (d)+(0.3,0) circle (0.03);
		\draw[fill] (d)+(-0.3,0) circle (0.03);
}}
\end{equation}

\paragraph{Remark.} This is the first case when an identity is valid exclusively under the action of the \( \Phi_L \) map. From the diagrammatic standpoint there is no way to multiply an edge by \( 2\hat{t} \), since the edge lacks any vertices.

\subsection{Tensor powers of the adjoint representation} 
\label{sec:adj-powers}

Vogel's framework is well suited for the adjoint sector, i.e. representations that are made of the tensor powers of the adjoint representation. In this case one should represent a Wilson line colored with the $ \text{adj}^{\otimes n} $ representation as $ n $ parallel lines. A vertex is replaced by the following sum (similar to the Leibnitz rule):

\begin{figure}[h!]
\centering
\mbox{\tikz[use Hobby shortcut, baseline = (current bounding box.center)]{
		\coordinate (a) at (0,0);
		\coordinate (b) at (4,0);
		\node at ($0.5*(a)+0.5*(b)+(0, 0.2)$) {\LARGE \(  = \sum \limits_{i=1}^{n}  \)  };
		\draw[-{Stealth[length=3mm]}, ultra thick] ($(a)+(-1,0)$) -- ($(a)+(1,0)$);
		\draw (a)-- +(0,1);
		\node at ($(a)+(-0.7, -0.3)$) {\(  \text{adj}^{\otimes n}  \) }; 
		\draw (b) ++ (-1, 0.25) -- ++ (2, 0);
		\draw (b) ++ (-1, 0.2) -- ++ (2, 0);
		\draw (b) ++ (-1, 0.15) -- ++ (2, 0);
		\draw (b) ++ (-1, -0.1) -- ++ (2, 0);
		\draw (b) ++ (-1, -0.15) -- ++ (2, 0);
		\draw (b) ++ (-1, -0.2) -- ++ (2, 0);
		\draw (b) ++ (-1, 0) -- ++ (2, 0);
		\draw (b)-- +(0,1);
		\draw[fill] (b) circle(0.03);
		\node at ($ (b)+(-1.1, 0) $) {\tiny \(  i  \) };
}}
\end{figure}

One may also consider a subspace corresponding to a Young diagram. In this case one needs to insert a projector, that is, up to a numeric factor, a Young symmetrizer acting on the strands of the diagram. 

For example, the second Casimir operator acting on the $ \text{adj}^{\otimes 2}$ can be diagrammatically expessed by \eqref{eq:2ndCas-on-square-1}:

\begin{equation}\label{eq:2ndCas-on-square-1}
\mbox{\tikz[scale=0.77, baseline = (current bounding box.center)]{
		\coordinate (a) at (0,0);
		\coordinate (b) at (3,0);
		\coordinate (c) at (6,0);
		\coordinate (d) at (9,0);
		\coordinate (e) at (12,0);
		\coordinate (f) at (15,0);
		\coordinate (g) at (18,0);
		\node at ($ 0.5*(a)+0.5*(b) $) {\(  =  \) };
		\node at ($ 0.5*(b)+0.5*(c) $) {\(  +  \) };
		\node at ($ 0.5*(c)+0.5*(d) $) {\(  +  \) };
		\node at ($ 0.5*(d)+0.5*(e) $) {\(  +  \) };
		\draw[-{Stealth[length=3mm]}, ultra thick] ($(a)+(-0.9,0)$) -- ($(a)+(0.9,0)$);
		\draw ($ (a) + (0:0.5)+(-0.15,0) $) arc(0:180:0.5);
		\node at ($(a)+(-0.7, -0.3)$) {\small \(   \text{adj}^{\otimes 2}   \) }; 
		\draw (b) +(-0.9, 0.2) -- +(0.9, 0.2);
		\draw (b) +(-0.9, -0.2) -- +(0.9, -0.2);
		\draw (b) +(0.24, 0.2) arc(0:180:0.24);
		\draw (c) +(-0.9, 0.2) -- +(0.9, 0.2);
		\draw (c) +(-0.9, -0.2) -- +(0.9, -0.2);
		\draw (c) +(-0.24, 0.2) .. controls +(0.05, 0.4) and +(0.0, 0.5) .. +(0.24, -0.2);
		\draw (d) +(-0.9, 0.2) -- +(0.9, 0.2);
		\draw (d) +(-0.9, -0.2) -- +(0.9, -0.2);
		\draw (d) +(+0.24, 0.2) .. controls +(-0.05, 0.4) and +(0, 0.5) .. +(-0.24, -0.2);
		\draw (e) +(-0.9, 0.2) -- +(0.9, 0.2);
		\draw (e) +(-0.9, -0.2) -- +(0.9, -0.2);
		\draw (e) +(0.24, -0.2) arc(0:180:0.24);
}}
\end{equation}

Basically, this equality is the diagrammatic way to write the coproduct of the Casimir operator that corresponds to its action on the tensor square of the adjoint representation:

\begin{equation}\label{eq:coproduct_Cas}
\triangle(C_2) = \triangle(X_iX^i) = (X_i\otimes1 + 1\otimes X_i)(X^i\otimes1 + 1\otimes X^i) = 2X_i\otimes X^i + 1\otimes X_iX^i + X_iX^i\otimes 1
\end{equation}

Since $ X_i X^i$ in the adjoint representation is a multiplication by $ 2t$, $1\otimes X_iX^i + X_iX^i\otimes 1 = 4t \cdot 1\otimes 1$, diagrammatic equality  \eqref{eq:2ndCas-on-square-1} can be further simplfied, using the non-diagrammatic factorisation of 'bubble' as \(2 t\), valid under the action of \(\Phi_L \):

\begin{equation}\label{eq:2ndCas-on-square-2}
\mbox{\tikz[scale=0.77, baseline = (current bounding box.center)]{
		\coordinate (a) at (0,0);
		\coordinate (b) at (3,0);
		\coordinate (c) at (6,0);
		\coordinate (d) at (9,0);
		\coordinate (e) at (12,0);
		\coordinate (f) at (15,0);
		\coordinate (g) at (18,0);
		\node at ($ 0.5*(b)+0.5*(c) $) {\(  +  \) };
		\node at ($ 0.5*(c)+0.5*(d) $) {\(  +  \) };
		\node at ($ 0.5*(d)+0.5*(e) $) {\(  +  \) };
		\node at ($ 0.5*(e)+0.5*(f) $) {\(  = 4t \) };
		\node at ($ 0.5*(f)+0.5*(g) $) {\(  -2  \) };
		\draw (b) +(-0.9, 0.2) -- +(0.9, 0.2);
		\draw (b) +(-0.9, -0.2) -- +(0.9, -0.2);
		\draw (b) +(0.24, 0.2) arc(0:180:0.24);
		\draw (c) +(-0.9, 0.2) -- +(0.9, 0.2);
		\draw (c) +(-0.9, -0.2) -- +(0.9, -0.2);
		\draw (c) +(-0.24, 0.2) .. controls +(0.05, 0.4) and +(0.0, 0.5) .. +(0.24, -0.2);
		\draw (d) +(-0.9, 0.2) -- +(0.9, 0.2);
		\draw (d) +(-0.9, -0.2) -- +(0.9, -0.2);
		\draw (d) +(+0.24, 0.2) .. controls +(-0.05, 0.4) and +(0, 0.5) .. +(-0.24, -0.2);
		\draw (e) +(-0.9, 0.2) -- +(0.9, 0.2);
		\draw (e) +(-0.9, -0.2) -- +(0.9, -0.2);
		\draw (e) +(0.24, -0.2) arc(0:180:0.24);
		\draw (f) +(-0.9, 0.2) -- +(0.9, 0.2);
		\draw (f) +(-0.9, -0.2) -- +(0.9, -0.2);
		\draw (g) +(-0.9, 0.2) -- +(0.9, 0.2);
		\draw (g) +(-0.9, -0.2) -- +(0.9, -0.2);
		\draw (g) +(0, -0.2) -- +(0, 0.2);
}}
\end{equation}

\section{Kernel of Lie algebra weight system}
\label{sec:kernel}

We already discussed in Sections \ref{sec:char} and \ref{sec:weght-system} that simple Lie algebras induce weight systems on trivalent diagrams with $k$ external legs and/or with distinguished cycles:
\begin{equation}
\Phi_{L}: \Lambda \rightarrow L^{\otimes k}, \
\put(23,3){\line(0,1){13}}
\put(23,3){\line(3,-2){12.5}}
\put(23,3){\line(-3,-2){12.5}}
\put(23,3){\circle*{2}}
\hspace{12mm} \mapsto \ f^{abc},
\put(27,0){\line(0,1){15}}
\put(27,0){\circle*{2}}
\linethickness{0.5mm}
\put(12,0){\line(1,0){30}}
\hspace{18mm} \mapsto \ \tau(X^a),
\end{equation}
where $\tau$ is any finite-dimensional irreducible representation of $L$.

One might ask a question whether all the diagrammatic stuff carry any independent value or it is enough to know all the answers for all the simple Lie algebras. This question can be formalised in a following way:

\textbf{Does the Lie algebra weight system have a non-trivial kernel on some class of trivalent diagrams?}

The answer is positive. Basically, there are several ways to construct such diagrams in the kernel of all the simple Lie algebras. 

\begin{enumerate}
\item The first class of such diagrams are those that consist of the zero divisor $P_{15}$ together with a bubble $\bbl$ in the cases where the bubble is not factorisable  as an element in $\Lambda$:

	\begin{equation}\label{eq:kernel-examples}
\mbox{\tikz[scale=0.70, baseline = (current bounding box.center)]{
\coordinate (a) at (0,0);
\coordinate (b) at (5,0);
\coordinate (c) at (11,0);

	\draw ($(a)+(-156:1)$) arc(-156:160:1);
	\node (a1) at ($(a)+(-0.9,0)$) {\( \hat P_{15}\)};
	\draw (a1)--($(a)+(1,0)$);

	\node (b1) at ($(b)+(-0.5, 0)$) {\small \( \hat P_{15}\)};
		\draw (b1) -- +(-1,0);
		\draw ($(b)+(150:0.5)$) arc(150:-140:0.5);
		\draw (b)++(0.5,0) -- +(1, 0);

		\draw[ultra thick] (c) circle(1);
		\node (c2) at ($(c)+(-1.9,0)$) {\( D_{17}:= \) };
		\node (c1) at (c) {\( \hat P_{15}\)};
		\draw (c1) -- +(90:1);
		\draw (c1) -- +(210:1);
		\draw (c1) -- +(330:1);
	
}}
	\end{equation}
These diagrams are undetectable by the weight systems of the simple Lie (super)algebras. Although the bubble corresponds to $2t$, this is not a purely diagrammatic correspondence, valid only under the action of a weight system.

These diagrams appear, for example, in the Kontsevich integral. The latter diagram was constructed by Vogel in his seminal paper \cite{vogel2011algebraic} to prove that there exist elements of the chord diagram algebra that are not detected by weight systems of semisimple Lie superalgebras.

\item The second class are diagrams that have multiple connected components. In this case the weight system mixes the $\Lambda$-factors from different connected components all together:

	\begin{equation}\label{eq:kernel-examples-2}
\mbox{\tikz[scale=0.77, baseline = (current bounding box.center)]{
\coordinate (a) at (0,0);
\coordinate (b) at (4,0);
\node at ($0.5*(a)+0.5*(b)+(-0.1, -0.5)$) {\(- \)};

\draw[-{Stealth[length=3mm]}, ultra thick] (a)+(-1.5, -1) -- +(1.5,-1);
\node[inner sep=0.4mm] (a1) at ($(a)+(-0.8, 0)$) {\small \( \hat u \)};
\node[inner sep=0.4mm] (a2) at ($(a)+(0.5, 0)$) {\small \( \hat v \)};
\draw (a1) -- ($(a)+(-1.2, -1)$);
\draw (a1) -- ($(a)+(-0.8, -1)$);
\draw (a1) -- ($(a)+(-0.4, -1)$);
\draw ($(a2)+(0, -0.25)+(120:0.25)$) arc(120:420:0.25);
\draw (a2) -- +(0,-0.5);
\draw (a2)++(0, -0.25)+(-120:0.25) -- +(-0.3,-0.75); 
\draw (a2)++(0, -0.25)+(-60:0.25) -- +(0.3,-0.75);

\draw[-{Stealth[length=3mm]}, ultra thick] (b)+(-1.5, -1) -- +(1.5,-1);
\node[inner sep=0.4mm] (b1) at ($(b)+(-0.8, 0)$) {\small \( \hat v \)};
\node[inner sep=0.4mm] (b2) at ($(b)+(0.5, 0)$) {\small \( \hat u \)};
\draw (b1) -- ($(b)+(-1.2, -1)$);
\draw (b1) -- ($(b)+(-0.8, -1)$);
\draw (b1) -- ($(b)+(-0.4, -1)$);
\draw ($(b2)+(0, -0.25)+(120:0.25)$) arc(120:420:0.25);
\draw (b2) -- +(0,-0.5);
\draw (b2)++(0, -0.25)+(-120:0.25) -- +(-0.3,-0.75); 
\draw (b2)++(0, -0.25)+(-60:0.25) -- +(0.3,-0.75); 
}}
\end{equation}

\end{enumerate}

These diagrams also appear in the Kontsevich integral. However, since they are not primitive elements of the Jacobi diagram algebra, they do not pose a significant problem. The coefficients of these diagrams in the Kontsevich integral are Vassiliev invariants, which are also non‑primitive \cite{alvarez1997primitive}; that is, they can be expressed as products of Vassiliev invariants of lower orders. A discussion of non‑primitive diagrams can be found in Section 4.2.2 of \cite{khudoteplov2024construction}.

\bigskip

The presence of a kernel in the weight systems generated by semisimple Lie algebras is very important. It indicates that semisimple Lie algebras are \textbf{insufficient} to generate the entire space of functions on the diagrams under consideration.

\section{Universal formulae}
\label{sec:universal-formulae}

In this Section we show how to apply the diagrammatic technique introduced in the previous sections to calculate specific quantities. 

\subsection{Dimension}
\label{sec:dimension}

Lie algebra dimension is one of the first quantities that were found to admit universal expression \cite{vogel2011algebraic}. The key point in the derivation is the identity $ \dim_{\text{adj}}  = \tr \text{Id}_{\text{adj}} $. In terms of the diagrams trace corresponds to a loop:

\begin{equation}\label{eq:dim_as_circle}
\mbox{\tikz[baseline = (current bounding box.center)]{
		\draw (0,0) circle(0.5);
		\node at (-1.55, 0) {\(  \dim_{L}= \Phi_L \circ \) };
}}
\end{equation}

Now the goal is to find a diagrammatic formula that involves such a loop. Thankfully, such formula can be found in Vogel's paper \cite{vogel2011algebraic}:

\begin{equation}\label{eq:sigma_action_1}
\mbox{\tikz[scale=0.56, baseline = (current bounding box.center)]{
		\coordinate (a) at (-9, 0);
		\coordinate (b) at ($(a)+(3.25, 0)$);
		\coordinate (c) at ($(b)+(2, 0)$);
		\coordinate (d) at ($(c)+(4.5, 0)$);
		\coordinate (e) at ($(d)+(3.25, 0)$);
		\coordinate (f) at ($(e)+(2, 0)$);
		\node at ($(a)+(-2.25, 0)$) { \(\sigma \Phi_L \) };
		\node at ($(a)+(-1.1, 0)$) {\Huge \( ( \) };
		\node at ($0.55*(a)+0.45*(b)$) {\( -t \) };
		\node at ($(b)+(-0.95, 0)$) {\huge \( ( \, \) };
		\node at ($0.5*(b)+0.5*(c)$) {\( - \) };
		\node at ($(c)+(0.95, 0)$) {\huge \( \,) \) };
		\node at ($(c)+(1.1, 0)$) {\Huge \(  \, \ ) \  \) };
		\node at ($0.5*(c)+0.5*(d)$) {\( =\,\Phi_L \) };
		\node at ($(d)+(-1.1, 0)$) {\Huge \( ( \) };
		\node at ($0.55*(d)+0.45*(e)$) {\( -\omega \) };
		\node at ($(e)+(-0.95, 0)$) {\huge \( (\, \) };
		\node at ($0.5*(e)+0.5*(f)$) {\( - \) };
		\node at ($(f)+(0.95, 0)$) {\huge \( \,) \) };
		\node at ($(f)+(1.1, 0)$) {\Huge \( \,) \) };
		\draw (a)+(45:-1)--+(-0.4, 0)--+(135:1);
		\draw (a)+(-0.4, 0) -- +(0.4, 0);
		\draw (a)+(45:1) -- +(0.4, 0) -- +(135:-1);
		\draw (b)+ (135:1)-- +(45:1);
		\draw (b) +(135:-1)-- +(45:-1);
		\draw (c)+(45:1) -- +(45:-1);
		\draw (c)+(135:-1)-- +(135:-0.1);
		\draw (c) +(135:0.1) -- +(135:1);
		\draw (d)+(45:-1)--+(-0.4, 0)--+(135:1);
		\draw (d)+(-0.4, 0) -- +(0.4, 0);
		\draw (d)+(45:1) -- +(0.4, 0) -- +(135:-1);
		\draw [fill] (d) circle(0.1);
		\node[scale=0.7] at ($(d)+(0, -0.3)$) {\( 0 \) };
		\draw (e)+ (135:1)-- +(45:1);
		\draw (e) +(135:-1)-- +(45:-1);
		\draw (f)+(45:1) -- +(45:-1);
		\draw (f)+(135:-1)-- +(135:-0.1);
		\draw (f) +(135:0.1) -- +(135:1);
}}
\end{equation}

A few comments on the formula \eqref{eq:sigma_action_1} are needed. First, it contains a bullet point marked by \(0\), which is a placeholder for a specific linear sum of diagrams (see Appendix in \cite{khudoteplov2024construction}). Second, this is not really a diagrammatic equation as it requires some extra conditions, valid only under the action of \( \Phi_L \). These preconditions are listed in (Theorem 6.3 in \cite{vogel2011algebraic}). Third, equation \eqref{eq:sigma_action_1} does not involve a loop, but can be made so once we connect top two legs by an edge, which corresponds to taking a half-trace over one of the tensor factors:

\begin{equation}
\mbox{\tikz[scale=0.56, baseline = (current bounding box.center)]{
		\coordinate (a) at (-9, 0);
		\coordinate (b) at ($(a)+(3.25, 0)$);
		\coordinate (c) at ($(b)+(2, 0)$);
		\coordinate (d) at ($(c)+(4.5, 0)$);
		\coordinate (e) at ($(d)+(3.25, 0)$);
		\coordinate (f) at ($(e)+(2, 0)$);
		\node at ($(a)+(-2.25, 0)$) { $\sigma \Phi_L$};
		\node at ($(a)+(-1.1, 0)$) {\Huge $($};
		\node at ($0.55*(a)+0.45*(b)$) {$-t$};
		\node at ($(b)+(-0.95, 0)$) {\huge $(\,$};
		\node at ($0.5*(b)+0.5*(c)$) {$-$};
		\node at ($(c)+(0.95, 0)$) {\huge $\,)$};
		\node at ($(c)+(1.1, 0)$) {\Huge $ \, \ ) \ $};
		\node at ($0.5*(c)+0.5*(d)$) {$=\,\Phi_L$};
		\node at ($(d)+(-1.1, 0)$) {\Huge $($};
		\node at ($0.55*(d)+0.45*(e)$) {$-\omega$};
		\node at ($(e)+(-0.95, 0)$) {\huge $(\,$};
		\node at ($0.5*(e)+0.5*(f)$) {$-$};
		\node at ($(f)+(0.95, 0)$) {\huge $\,)$};
		\node at ($(f)+(1.1, 0)$) {\Huge $\,)$};
		\draw (a)+(45:-1)--+(-0.4, 0);
		\draw (a)+(-0.4, 0) -- +(0.4, 0);
		\draw (a)+(0.4, 0) -- +(135:-1);
		\draw ($ (a)+(0.4, 0) $) arc(-30:210:0.465);
		\draw ($ (b)+(0,0.15) $) circle(0.45);
		\draw (b) +(135:-1)-- +(45:-1);
		\draw (c)+(135:-1)-- +(135:-0.1);
		\draw (c) +(135:0.1) -- +(135:0.3);
		\draw (c)+(45:0.3) -- +(45:-1);
		\draw ($ (c)+(45:0.3) $) arc(-45:225:0.3 );
		\draw (d)+(45:-1)--+(-0.4, 0);
		\draw (d)+(-0.4, 0) -- +(0.4, 0);
		\draw (d)+(0.4, 0) -- +(135:-1);
		\draw [fill] (d) circle(0.1);
		\node[scale=0.7] at ($(d)+(0, -0.3)$) {$0$};
		\draw ($ (d)+(0.4, 0) $) arc(-30:210:0.465);
		\draw ($ (e)+(0,0.15) $) circle(0.45);
		\draw (e) +(135:-1)-- +(45:-1);
		\draw (f)+(135:-1)-- +(135:-0.1);
		\draw (f) +(135:0.1) -- +(135:0.3);
		\draw (f)+(45:0.3) -- +(45:-1);
		\draw ($(f)+(45:0.3)$) arc(-45:225:0.3);
}}
\end{equation}

The bubble in the first term can be factorised as $ 2 t $ and the diagram with a bullet marked by zero stands for a linear combination of diagrams that turns out to be zero: 

\begin{equation}
	\mbox{\tikz[scale=0.56, baseline = (current bounding box.center)]{
			\coordinate (a) at (4,0);
			\coordinate (b) at ($(a)+(3.5,0)$);
			\coordinate (c) at ($(b)+(3.2,0)$);
			\coordinate (d) at ($(c)+(3.1,0)$);
			\coordinate (e) at ($(d)+(3.1,0)$);
			
			\node at ($(b)+(-1.5,.07)$) {\( =\;2t^2 \)};
			\node at ($(c)+(-1.5,.07)$) {\( +\;2t\)};
			\node at ($(d)+(-1.4,.07)$) {\( -\; \frac43\)};
			\node at ($(e)+(-1.4,.07)$) {\( -\; \frac23\)};
			\node at ($(e)+(1.7,0.07)$) {\( =\;0\)};
			
			\draw (b)+(45:-1)--+(-0.4, 0);
			\draw (b)+(-0.4, 0) -- +(0.4, 0);
			\draw (b)+(0.4, 0) -- +(135:-1);
			\draw ($ (b)+(0.4, 0) $) arc(-30:210:0.465);
			
			\draw (a)+(45:-1)--+(-0.4, 0);
			\draw (a)+(-0.4, 0) -- +(0.4, 0);
			\draw (a)+(0.4, 0) -- +(135:-1);
			\draw [fill] (a) circle(0.1);
			\node[scale=0.7] at ($(a)+(0, -0.3)$) {$0$};
			\draw ($ (a)+(0.4, 0) $) arc(-30:210:0.465);
			
			\draw (c)+(45:-1)--+(-0.4, 0);
			\draw (c)+(-0.4, 0) -- +(0.4, 0);
			\draw (c)+(0.4, 0) -- +(135:-1);
			\draw (c)++(0.4, 0)+(112:-0.3) -- ($(c)+(-0.4, 0)+(68:-0.3)$);
			\draw ($ (c)+(0.4, 0) $) arc(-30:210:0.465);
			
			\draw (d)+(45:-1)--+(-0.4, 0);
			\draw (d)+(-0.4, 0) -- +(0.4, 0);
			\draw (d)+(0.4, 0) -- +(135:-1);
			\draw (d)++(0.4, 0)+(112:-0.4) -- ($(d)+(-0.4, 0)+(68:-0.4)$);
			\draw (d)++(0.4, 0)+(112:-0.2) -- ($(d)+(-0.4, 0)+(68:-0.2)$);
			\draw ($ (d)+(0.4, 0) $) arc(-30:210:0.465);
			
			\draw (e)+(45:-1)--+(-0.4, 0);
			\draw (e)+(-0.4, 0) -- +(0.4, 0);
			\draw (e)+(0.4, 0) -- +(135:-1);
			\draw (e)++(0.4, 0)+(112:-0.3) -- ($(e)+(-0.4, 0)+(68:-0.3)$);
			\draw (e) -- +(0, -0.27);
			\draw ($ (e)+(0.4, 0) $) arc(-30:210:0.465);
			
	}}
\end{equation}

Hence, the whole equation can be rewritten as:

\begin{equation}
\sigma (2t - t(\text{dim}_L-1))=-\omega((\text{dim}_L-1))
\end{equation} 

It follows that

\begin{equation}\label{eq:univ-dim}
\text{dim}_L = \frac{\omega -3 t \sigma}{\omega - t \sigma}=\frac{(\alpha-2t)(\beta-2t)(\gamma-2t)}{\alpha \beta \gamma}
\end{equation}

\subsection{Casimir invariants}
\label{sec:Casimirs}

Another example of a universal formula is an expression for the higher Casimir eigenvalues in adjoint representation that was found in paper \cite{mkrtchyan2012casimir}. Analogous result was achieved for Lie superalgebras in \cite{isaev2022split}. Both papers refer to the Okubo's paper \cite{okubo1977casimir}, which contains a relation between higher Casimir eigenvalues and the second Casimir eigenvalues on irreducible representations in the adjoint tensor square.

The higher Casimir invariants in these papers are defined by:

\begin{equation}\label{eq:higher_Cas_def}
C_p = \Tr_{ad}(X^{\mu_1}X^{\mu_2}\cdots X^{\mu_p}) X_{\mu_1}X_{\mu_2}\cdots X_{\mu_p}
\end{equation}

The claim \cite{mkrtchyan2012casimir, isaev2022split} is that with $t_2 = \alpha^2+\beta^2+\gamma^2$ and $t_3 = \alpha^3+\beta^3+\gamma^3$ a generating function $C(z) = \sum_{p=0}^\infty C_p z^p$ becomes:
\begin{multline}
\label{eq:gen_fun}
C(z) = \frac{(\alpha-2t)(\beta-2t)(\gamma-2t)}{\alpha \beta \gamma} +z^2 \frac{96t^3+168t^3 z+6(14t^3+t t_2-t_3)z^2+(13t^3+3 t t_2-4t_3)z^3}{6(2t+\alpha z)(2t+\beta z)(2t+\gamma z)(2+z)(1+z)}
\end{multline}

This result becomes very apparent once the diagrammatic technique is used. Indeed, an operator $ C_p $ as in (\ref{eq:higher_Cas_def}) can be represented by a following diagram:

\begin{equation}\label{eq:higher_Cas_as_diag}
\mbox{\tikz[baseline = (current bounding box.center)]{
		\node at (-1.5, 0.75) {$(-1)^p C_p=$};
		\draw (-1, 0) -- (1, 0);
		\draw (0, 1) circle (0.5);
		\draw ($(0, 1)+(50:-0.5)$) -- (-0.5, 0);
		\draw ($(0, 1)+(60:-0.5)$) -- (-0.4, 0);
		\draw ($(0, 1)+(70:-0.5)$) -- (-0.3, 0);
		\draw ($(0, 1)+(80:-0.5)$) -- (-0.2, 0);
		\node at (0.1, 0.25) {$\cdot$};
		\node at (0., 0.25) {$\cdot$};
		\node at (0.2, 0.25) {$\cdot$};
		\draw ($(0, 1)+(120:-0.5)$) -- (0.4, 0);
		\draw ($(0, 1)+(130:-0.5)$) -- (0.5, 0);
		\node at (1, 0.75) {$=$};
		\node (a) at (2.5, 0.75) {$\hat{x}_{p-1}$};
		\draw (1.5, 0.75) -- (a);
		\draw ($(3, 0.75)+(160:0.75)$) arc(160:-160:0.75);
		\draw (3.75, 0.75) -- (4.5, 0.75);
		\node at (4.75, 0.75) {$=$};
		\node at (5.4, 0.75) {$\hat{x}_{p-1} \cdot$};
		\draw (6, 0.75) -- +(0.5, 0);
		\draw (7, 0.75) circle (0.5);
		\draw (7.5, 0.75) -- (8, 0.75);
		\node at (9, 0.75) {$=2 t \hat{x}_{p-1} \cdot$};
	\draw (10, 0.75) -- +(2, 0);
}}
\end{equation}

From the equation (\ref{eq:higher_Cas_as_diag}) it is clear that there is a universal formula for $ C_p $:

\begin{equation}
C_p = 2\, (-1)^p \, t\, \chi_L(\hat{x}_{p-1})
\end{equation}

This expression is very compact, but formulated in terms of an infinitely many elements of $ \Lambda $ algebra $\{\hat{x}_n\} $.

To express $ \hat{x}_n $ in terms of $ \alpha $, $ \beta $ and $ \gamma $ one can use a relation discovered by Jan Kneissler \cite{kneissler2001spaces}, which is given in a slightly reformulated form in paper \cite{vogel1999universal}:
\begin{multline}\label{eq:kneiss-relation}
\chi_L({\hat{x}_{n+3}}) = (\alpha + \beta + \gamma) \chi_L({\hat{x}_{n+2}}) - (\alpha \beta + \beta \gamma + \gamma \alpha) \chi_L({\hat{x}_{n+1}}) + \alpha \beta \gamma \chi_L({\hat{x}_{n}}) +\\+ \frac{(\alpha \beta + \beta \gamma + \gamma \alpha) (\alpha + \beta + \gamma)^{n+1}}{2}-\frac{\alpha \beta \gamma (\alpha + \beta + \gamma)^n}{2} - \alpha \beta \gamma (2(\alpha + \beta + \gamma))^n
\end{multline}

Before applying (\ref{eq:kneiss-relation}), one needs to remember that $\hat{x}_0=0$, $\hat{x}_1=2\hat{t}$ and $\hat{x}_2 = \hat{t}^2$. Next, one can solve an equation (\ref{eq:kneiss-relation}) as a linear recurrence relation. In this case the solution is:

\begin{equation}
\chi_L({\hat{x}_n}) = -\frac{\alpha \beta \gamma(2t)^n}{(2t-\alpha)(2t-\beta)(2t-\gamma)}+\frac{t^n}{2} + c_\alpha \alpha^n + c_\beta \beta^n + c_\gamma \gamma^n
\end{equation}

The values $c_\alpha, c_\beta, c_\gamma$ are determined by $\hat{x}_0=0$, $\hat{x}_1=2\hat{t}$ and $\hat{x}_2 = \hat{t}^2$ conditions:

\begin{equation}
c_\alpha=-\frac{-t^2(\dim-8)+(2+\dim)\beta \gamma +(3 \dim - 4)t(\beta+\gamma)}{2\dim (\alpha-\beta)(\alpha-\gamma)}
\label{eq:c-alpha}	
\end{equation}

Here $\dim = \frac{(\alpha-2t)(\beta-2t)(\gamma-2t)}{\alpha \beta \gamma}$ as in formula \eqref{eq:univ-dim} for the universal Lie algebra dimension. Expression for $c_\beta$ can be derived from the formula (\ref{eq:c-alpha}) for the $c_\alpha$ by a replacement $\alpha \leftrightarrow \beta$, and the formula for $c_\gamma$ by the $\alpha \leftrightarrow \gamma$ replacement.

Comparison with (\ref{eq:gen_fun}) shows that

\begin{equation}
C(z) = \dim + \sum_{k=1}^{\infty} \frac{\chi_L(\hat{x}_k)}{(2t)^k} (-z)^{k+1}
\end{equation}

It means that Mkrtchyan's result for the $ C_p $ deviates by a factor of $ (2t)^p $, which indicates that in the paper \cite{mkrtchyan2012casimir} the Killing form was used as a metric so that $ t=1/2 $.

\subsection{Decomposition of the adjoint tensor powers}
\label{sec:adj-square}

In this subsection we discuss the representation theory for the adjoint sector. There have been many developments in this field recently, going up to the 5th tensor power \cite{avetisyan2024uniform, isaev2024split}. It also was the primary focus of the Vogel's 1999 preprint \cite{vogel1999universal}. Despite all the interest towards this topic, it remains quite vague from the diagrammatic viewpoint.

Remember the example of the second Casimir operator acting on the tensor square of the adjoint representation \eqref{eq:2ndCas-on-square-1} from the Section \ref{sec:adj-powers}. We have seen that it could be split into two parts. One of them is the 'bubble' insertion which factors out as $ 2 t $ under \(\Phi_L\), and the other one is non-trivial and represented by a diagram  $ \Psi = {\text{\Large $\sqsupset \hspace{-0.95mm} \sqsubset$} }$. 

Basically, this diagram correspond to an endomorphism of $ \text{adj}^{\otimes 2} $ given by $ x \otimes y \mapsto \sum_i [x, e_i] \otimes [e^i, y] $ . This operator is called a split Casimir operator in literature \cite{isaev2022split} (in these papers the split Casimir operator differs by a sign from our $ \Psi = {\text{\Large $\sqsupset \hspace{-0.95mm} \sqsubset$} }$).

The paradigm is to find a polynomial relation the operator satisfies, which would determine the eigenvalues of the split Casimir. Consequently, it would provide a decomposition into the eigenspaces. The dimensions of these eigenspaces can be calculated by taking the trace over the corresponding projectors.

\paragraph{Antisymmetric square} $\bigwedge^2(\text{adj})$ can be obtained by the following projector:
\begin{equation}
\mbox{\begin{picture}(-20,10)
		\put(-100,0){$P_a = \dfrac{1}{2} \Big($}
		\put(-60,-3){\line(1,0){20}}
		\put(-60,9){\line(1,0){20}}
		\put(-35,0){$-$}
		\put(-20,-3){\line(5,2){13}}
		\put(10,9){\line(-5,-2){13}}
		\put(-20,9){\line(5,-2){30}}
		\put(15,0){$\Big),$}
		\put(35,0){$P_a^2 = P_a$}
\end{picture}}
\end{equation}

$\bigwedge^2(\text{adj})$ is the easiest but illustrative case. The restriction of $ \Psi = {\text{\Large $\sqsupset \hspace{-0.95mm} \sqsubset$} }$ onto the antisymmetric square takes the form $P_a \Psi P_a = \Psi P_a =  \frac12 {\text{\Large $\sqsupset \hspace{-0.95mm} \sqsubset$} } - \frac12 {\text{\Large $\sqsupset   \hspace{-0.95mm} > \hspace*{-0.99 mm}< $}} = \frac12 {\text{\Large$>\hspace*{-2.5mm}-\hspace*{-2.5mm}<$}}$, where the IHX relation has been applied.

Operator \text{\Large $ {>\hspace*{-2.5mm}-\hspace*{-2.5mm}<} $} satisfies a polynomial relation $ {\text{\Large$ >\hspace*{-2.5mm}-\hspace*{-2.5mm}<$}}^2 - 2 t {\text{\Large $>\hspace*{-2.5mm}-\hspace*{-2.5mm}<$}} =0 $, since $ {\text{\Large $>\hspace*{-2.5mm}-\hspace*{-2.5mm}<$}}^2 = {\text{\Large $>\hspace*{-2.5mm}-\hspace*{-2.5mm}<$}}\hspace{-1.0mm}{\text{\Large$>\hspace*{-2.5mm}-\hspace*{-2.5mm}<$}} $ and the 'bubble' factorizes as $ 2 t $ (this time it is a purely diagrammatic factorisation).

At this point there is a fork depending on the value of \(t \). If \( t = 0\), as is the case with \(D(2,1,\lambda) \) superalgebra, \( \Psi \) is a nilpotent operator and \( \bigwedge^2(\text{adj})\) is reducible but indecomposable.  When \( t \neq 0 \), which is the case for all simple Lie algebras, the space $ \wedge^2 \text{adj} $ splits into the two eigenspaces $X_1$ and $X_2$, corresponding to eigenvalues  $ t $ and \( 0\) of $ \Psi = {\text{\Large $\sqsupset \hspace{-0.95mm} \sqsubset$} }$:

\begin{equation}
    \mbox{\begin{picture}(-20,10)
    \put(-94,0){$P_{\text{adj}} = $}
    \put(-60,9){\line(1,0){10}}
    \put(-60,-3){\line(1,0){10}}
    \put(-50,-3){\line(1,1){6}}
    \put(-50,9){\line(1,-1){6}}
    \put(-44,3){\line(1,0){10}}
    \put(-34,3){\line(1,1){6}}
    \put(-34,3){\line(1,-1){6}}
    \put(-28,9){\line(1,0){10}}
    \put(-28,-3){\line(1,0){10}}
    \put(-15,0){$,  \quad  P_{\text{adj}}^2 = 2\hat{t} \cdot P_{\text{adj}}$}
		    \end{picture}}
\end{equation}

\begin{equation}
	\bigwedge\nolimits^2 L = X_1 \oplus X_2, \qquad 
\end{equation}

The dimensions for the these spaces can be optained by taking trace of  \( P_a \Psi P_a = \frac12 {\text{\Large$>\hspace*{-2.5mm}-\hspace*{-2.5mm}<$}}\) yielding \( \frac12  \mbox{\tikz[scale=0.2, baseline = -0.8mm]{\draw (0,0) circle(1); \draw (-1,0) -- (1,0);}}\), which is mapped to \( t\, \text{dim}_L \) by \( \Phi_L\). Since \( \tr (P_a \Psi P_a) = t \cdot \text{dim}_{X_1} + 0\cdot  \text{dim}_{X_2} \), we may conclude that \( \text{dim}_{X_1} = \text{dim}_L \) and \( \text{dim}_{X_2} = \frac{\text{dim}_L (\text{dim}_L-3)}{2} \).

\paragraph{Symmetric square} $S^2(\text{adj})$ can be obtained by the following projector:
\begin{equation}
\mbox{\begin{picture}(-20,10)
		\put(-100,0){$P_s = \dfrac{1}{2} \Big( $}
		\put(-60,-3){\line(1,0){20}}
		\put(-60,9){\line(1,0){20}}
		\put(-35,0){$+$}
		\put(-20,-3){\line(5,2){13}}
		\put(10,9){\line(-5,-2){13}}
		\put(-20,9){\line(5,-2){30}}
		\put(15,0){$\Big),$}
		\put(35,0){$P_s^2 = P_s$}
\end{picture}}
\end{equation}

To begin with, remember that there exists a $ 1 $-dimensional subspace generated by the Casimir element $ e_i \otimes e^i $. It corresponds to a 1-dimensional representation (singlet) \( X_0 \) in the \( S^2 \text{adj} \). A mapping onto the \( X_0 \) can be represented by a following diagram:
\begin{equation}
\mbox{\begin{picture}(30,10)
		\put(-50,0){$P_{0} = $}
		\put(-20,3){\oval(14,14)[r]}
		\put(5,3){\oval(14,14)[l]}
		\put(8,0){$,$}
		\put(20,0){$P_{0}^2 = \bigcirc \cdot P_{0}$}
\end{picture}}
\end{equation}

The space $ X_0 $ is also an eigenspace for $ {\text{\Large $\sqsupset \hspace{-0.95mm} \sqsubset$} } $, since $ {\text{\Large $\sqsupset \hspace{-0.95mm} \sqsubset$} }\hspace{-2.5mm} \text{\Large$\supset \subset$} = 2t  \text{\Large$\supset \subset$} $, if we use \( \Phi_L \) to factor out the 'bubble'.

The leftover subspace is denoted by $ Y = S^2 (\text{adj}) / X_0 $. Unlike the previous cases, $ Y $ is not an eigenspace for \( \Psi \), suggesting that it can be futher decomposed. Indeed, under a few additional assumptions, such as \( \Psi \) acting bijectively on \( Y \), it can be shown that \( \varphi = \Psi_L|_Y \) satisfyes a cubic relation on \( Y \):

\begin{equation}
\varphi^3 - t \varphi^2 + (\sigma - 2 t^2) \varphi - (\omega - t \sigma) =0
\end{equation}

Variables \( \alpha \), \( \beta \) and \( \gamma \) are introduced as roots of this polynomial, i.e. the eigenvalues of \( \varphi \), corresponding to eigenspaces \(Y_\alpha \), \( Y_\beta \) and \( Y_\gamma\), respectively.

\paragraph{Antisymmetric cube} Consider the following operators \( O_a \)  and \( O_s \):
\begin{equation}
\mbox{\tikz[scale=0.9, baseline = (current bounding box.center)]{
		\draw[rounded corners=8pt] (0.2,0) rectangle (1,2);
		\draw (0.2, 1) -- +(-0.8,0);
		\draw (0.2, 1.6) -- +(-0.8,0);
		\draw (0.2, 0.4) -- +(-0.8,0);
		\draw (1, 1.6) arc(90:-90:0.3);
		\draw (1.3, 1.3) -- (2.3, 1.3);
		\draw (1, 0.4) -- (2.3, 0.4);
		\draw[rounded corners=8pt] (2.3, 0) rectangle (3, 1.7);
		\draw (4, 1.3) -- (3, 1.3);
		\draw (4.3, 0.4) -- (3, 0.4);
		\draw (4.3, 1.6) arc (90:270:0.3);
		\draw[rounded corners=8pt] (4.3,0) rectangle (5.1, 2);
		\draw (5.1, 1) -- +(0.8,0);
		\draw (5.1, 1.6) -- +(0.8,0);
		\draw (5.1, 0.4) -- +(0.8,0);
		\node at (0.6,1) {\( P_a \)};
		\node at (2.65,0.85) {\footnotesize \( P_{a/s} \)};
		\node at (4.7,1) {\( P_a \)};
		\node at (-1.4, 1) {\( O_{a/s} = \)};
}}
\end{equation}

Their form suggests that they correspond to mappings from \( \wedge^3 \text{adj}\) to the subrepresentations appearing in \( \wedge^2 \text{adj}\) and \( S^2 \text{adj}\), respectively. If we allow for factorisation of the 'bubble' as \(2t\) these operators can be shown to satisfy \( O_a O_s = O_s O_a = 0\).

However, not all of the representations from \( \text{adj}^{\otimes 2}\) appear in \( \wedge^3\text{adj}\). Indeed, consider the \(X_1 \) representation in \( \wedge^2 \text{adj}\). As was previously described, it is possible to map onto \(X_1 \) by $ {\text{\large$>\hspace*{-2.2mm}-\hspace*{-2.2mm}<$}} $ operator. Next, if we insert $ {\text{\large$>\hspace*{-2.2mm}-\hspace*{-2.2mm}<$}} $ inside the \( O_a \) diagram, we would end up with zero diagram because of the IHX relation:

\begin{equation}
\mbox{\tikz[scale=0.9, baseline = (current bounding box.center)]{
		\draw[rounded corners=8pt] (0.2,0) rectangle (1,2);
		\draw (0.2, 1) -- +(-0.8,0);
		\draw (0.2, 1.6) -- +(-0.8,0);
		\draw (0.2, 0.4) -- +(-0.8,0);
		\draw (1, 1.6) arc(90:-90:0.3);
		\draw (4.3, 1.6) arc (90:270:0.3);
		\draw (1, 0.4) -- +(0.3,0);
		\draw (4.3, 0.4) -- +(-0.3,0);
		\draw (1.3, 1.3) arc(90:-90:0.45);
		\draw (4, 1.3) arc(90:270:0.45);
		\draw (1.75, 0.85) -- (3.55, 0.85);
		\draw[rounded corners=8pt] (4.3,0) rectangle (5.1, 2);
		\draw (5.1, 1) -- +(0.8,0);
		\draw (5.1, 1.6) -- +(0.8,0);
		\draw (5.1, 0.4) -- +(0.8,0);
		\node at (0.6,1) {\( P_a \)};
		\node at (4.7,1) {\( P_a \)};
		\node at (6.5,1) {\( =0\)};
}}
\end{equation}

Therefore, \( X_1 \) representation is not present in the image of the \( O_a \) operator. However, it does not guarantee that there is no other operator that maps \( \wedge^3 \text{adj} \) onto the \( X_1 \). For all simple Lie algebra this is the case, however, from the diagrammatic standpoint it depends on whether there exists a non-zero diagram of the form:

\begin{equation}
\mbox{\tikz[scale=0.9, baseline = (current bounding box.center)]{
		\draw[rounded corners=8pt] (0.2,0) rectangle (1,2);
		\draw (0.2, 1) -- +(-0.8,0);
		\draw (0.2, 1.6) -- +(-0.8,0);
		\draw (0.2, 0.4) -- +(-0.8,0);
		\draw (1., 1) -- +(0.4,0);
		\draw (1., 1.6) -- +(0.4,0);
		\draw (1., 0.4) -- +(0.4,0);
		\draw (1.4, 0) rectangle (3, 2);
		\draw (3, 1) -- (4.5, 1);
		\draw (4.5, 0) rectangle (6.1, 2);
		\draw (6.1, 0.4) -- +(0.4,0);
		\draw (6.1, 1.6) -- +(0.4,0);
		\draw (6.1, 1) -- +(0.4,0);
		\draw[rounded corners=8pt] (6.5,0) rectangle (7.3, 2);
		\draw (7.3, 1) -- +(0.8,0);
		\draw (7.3, 1.6) -- +(0.8,0);
		\draw (7.3, 0.4) -- +(0.8,0);
		\node at (0.6,1) {\( P_a \)};
		\node at (6.9,1) {\( P_a \)};
		\node at (2.2, 1) {\footnotesize subgraph};
		\node at (5.3, 1) {\footnotesize subgraph};
}}
\end{equation}

What makes it tricky is that the subgraph might belong to the kernel of all Lie algebra weight systems, as was the case in Section \ref{sec:kernel}. If so, then our analysis of the representation theory of the simple Lie algebras is irrelevant, since some of the subrepresentations may appear exclusively in the diagrammatic sense.

This example is a good illustration of the obstacles that could appear on the path towards a purely diagrammatic formulation of the universal decomposition. These obstacles are our lack of understanding of the spaces of diagrams with more than \(3 \) legs as well as the existence of diagrams, undetected by \( \Phi_L \).

Despite the aforementioned troubles, the diagrammatic technique is still quite useful for the universal decomposition, as we tried to demonstrate in this subsection.

\subsection{Quantum dimension and symmetric Casimir invariants}
\label{sec:qdim-and-symm-Cas}

\subsubsection{Quantum dimension}
\label{sec:qdim}

Quantum dimension is a quantum invariant of the unknot. This invariant can be calculated with the help of quantum R-matrix via the Reshetikhin-Turaev approach \cite{Reshetikhin:1990pr, reshetikhin1991invariants}. However it is unknown the Vogel's universalization for this approach (see attempts in \cite{mironov2016universalA}). On the other hand, any quantum knot invariant can be calculated by the applying of the corresponding weight system $\Phi$ to the Kontsevich integral:

\begin{equation}
\text{QI}_{L}^{R}(\mathcal{K}) = \Phi_{L}^{R} \Big( \, \text{I}(\mathcal{K}) \, \Big) 
\end{equation}

The Kontsevich integral for the unknot was calculated in \cite{bar2000wheels, bar2003two}:

\begin{equation}\label{eq:ZunknB}
\text{I}\left( \bigcirc \right) = \rho \circ \exp\left( \sum_{n=1}^{\infty} b_{2n}\, w_{2n} \right),
\end{equation}
where the 'wheels' $w_{2n}$ is the open Jacobi diagrams with $2n$ legs

\begin{equation}
\put(-175,0){$w_{2} = $}
\put(-130,3){\circle{20}}
\put(-150,3){\line(1,0){10}}
\put(-120,3){\line(1,0){10}}
\put(-105,0){$, \ \ w_{4} = $}
\put(-60,-10){\line(1,0){20}}
\put(-60,-10){\line(0,1){20}}
\put(-60,10){\line(1,0){20}}
\put(-40,-10){\line(0,1){20}}
\put(-60,-10){\line(-1,-1){7}}
\put(-60,10){\line(-1,1){7}}
\put(-40,-10){\line(1,-1){7}}
\put(-40,10){\line(1,1){7}}
\put(-25,0){$, \ \ w_{6} = $}
\put(20,-5){\line(0,1){15}}
\put(20,10){\line(2,1){13}}
\put(32.8,16.5){\line(2,-1){13}}
\put(45.8,-5){\line(0,1){15}}
\put(20,-5){\line(2,-1){13}}
\put(32.8,-11.5){\line(2,1){13}}
\put(20,-5){\line(-3,-2){8}}
\put(20,10){\line(-3,2){8}}
\put(32.8,16.5){\line(0,1){9}}
\put(32.8,-11.5){\line(0,-1){9}}
\put(45.8,-5){\line(3,-2){8}}
\put(45.8,10){\line(3,2){8}}
\put(60,0){$, \ \ w_{8} = $}
\put(110,-3){\line(0,1){10}}
\put(110,7){\line(1,1){7}}
\put(117,14){\line(1,0){10}}
\put(127,14){\line(1,-1){7}}
\put(134,7){\line(0,-1){10}}
\put(134,-3){\line(-1,-1){7}}
\put(110,-3){\line(1,-1){7}}
\put(117,-10){\line(1,0){10}}
\put(110,-3){\line(-2,-1){7}}
\put(110,7){\line(-2,1){7}}
\put(117,14){\line(-1,2){3.3}}
\put(127,14){\line(1,2){3.3}}
\put(134,7){\line(2,1){7}}
\put(134,-3){\line(2,-1){7}}
\put(117,-10){\line(-1,-2){3.3}}
\put(127,-10){\line(1,-2){3.3}}
\put(148,0){$, \ \ \ldots $}
\end{equation}

\vspace{4mm}

\noindent and the modified Bernoulli numbers $b_{2n}$ are defined by the power series expansion
\begin{gather}
\sum_{n=0}^{\infty} b_{2n}\, x^{2n} = \dfrac{1}{2}\log \dfrac{\sinh\frac{x}{2}}{\frac{x}{2}},  \\
b_2 = \dfrac{1}{48}, \ b_4 = -\dfrac{1}{5760}, \ b_6 = \dfrac{1}{362880}, \ b_8 = -\dfrac{1}{19353600}, \ \ldots \nn
\end{gather}

Diagrams $w_{2n}$ are open Jacobi diagram, while the Kontsevich integral is over chord diagrams or closed Jacobi diagrams. 
To represent it in  the corresponding form we use the isomorphism between two algebras, which is given by the symmetrization map:
\begin{align}
\rho: \mathcal{B} \rightarrow \mathcal{C}.
\end{align}

This map is defined as follows:
\begin{align}
\mbox{\begin{picture}(300,25)(-110,0)
		\put(-35,0){$\rho\Big($}
		\put(24,0){$\Big)$}
		\put(0,3){\circle{20}}
		\put(0,3){\color{lightgray}\circle*{19.8}}
		\put(10,3){\line(1,0){10}}
		\put(-20,3){\line(1,0){10}}
		\put(0,13){\line(0,1){10}}
		\put(0,-7){\line(0,-1){10}}
		\put(-7,10){\line(-1,1){7}}
		\put(-7,-4){\line(-1,-1){7}}
		\put(7,10){\line(1,1){7}}
		\put(6,-12){\circle*{1.5}}
		\put(11.5,-8){\circle*{1.5}}
		\put(15,-3){\circle*{1.5}}
		\put(15,-12){\tiny $n$}
		\put(33,0){$= \dfrac{1}{n!} \sum\limits_{\sigma \in S_{n}}$}
		\put(102,8){\oval(25,12)}
		\put(95.4,8){\color{lightgray}\circle*{11.5}}
		\put(108.6,8){\color{lightgray}\circle*{11.5}}
		\put(93,2.5){\line(0,-1){8.5}}
		\put(96,2){\line(0,-1){8}}
		\put(108,2){\line(0,-1){8}}
		\put(111,2.5){\line(0,-1){8.5}}
		\put(99,-2){\circle*{1}}
		\put(102,-2){\circle*{1}}
		\put(105,-2){\circle*{1}}
		\put(89.5,-6){\line(1,0){25}}
		\put(89.5,-14){\line(1,0){25}}
		\put(89.5,-6){\line(0,-1){8}}
		\put(114.5,-6){\line(0,-1){8}}
		\put(100,-12){\footnotesize $\sigma$}
		\put(92,-14){\line(0,-1){9}}
		\put(95,-14){\line(0,-1){10}}
		\put(112,-14){\line(0,-1){9}}
		\put(109,-14){\line(0,-1){10}}
		\put(99,-19){\circle*{1}}
		\put(102,-19){\circle*{1}}
		\put(105,-19){\circle*{1}}
		\linethickness{4.06mm}
		\put(94.5,8.01){\color{lightgray}\line(1,0){15}}
		\linethickness{0.4mm}
		\put(102,0){\circle{50}}
		\put(126.8,3){\vector(0,1){2}}
\end{picture}}
\end{align}

\vspace{8mm}

\noindent where a gray part contains an arbitrary trivalent graph. Thus we can get the Kontsevich integral for the unknot in terms of closed Jacobi diagrams:
\begin{align}
\label{eq:symmunk}
\text{I}\left( \bigcirc \right) = 
1 + b_2 \rho(w_2) + \left( \dfrac{b_2^2}{2!} \rho(w_2^2) + b_4 \rho(w_4) \right) + \left( \dfrac{b_2^3}{3!} \rho(w_2^3) + b_2b_4 \rho(w_2w_4) + b_6 \chi(w_6) \right)  + ...
\end{align}

The product of open Jacobi diagrams is their disjoint union.  Note that $\rho(w_nw_m) \neq  \rho(w_n) \rho(w_m)$.

Next step is to consider the Lie algebra weight system $\Phi_{L}^{R}$ in the adjoint representation $R=\text{adj}$
\begin{align}
\Phi_{L}^{\text{adj}} \Big( \text{I}\left( \bigcirc \right) \Big).
\end{align}

As we discussed in Section \ref{sec:rules} in order to consider the adjoint representation we substitute bold lines by ordinary lines (no distinguished cycle) and forget about the orientation. Therefore, the Kontsevich integral for the unknot in the adjoint representation is a sum over 3-graphs. For example, first few terms have the following form:

\begin{align}
\Phi_{L}^{\text{adj}} \Big( \text{I}\left( \bigcirc \right) \Big) \ = \ 
\dim L \ + \ \dfrac{b_2}{2!}\sum\limits_{\sigma \in S_2} 
\mbox{\put(22,3){\circle{40}}
	\put(22,9){\oval(20,14)}
	\put(18,2){\line(0,-1){5}}
	\put(26,2){\line(0,-1){5}}
	\put(15,-3){\line(1,0){14}}
	\put(15,-9){\line(1,0){14}}
	\put(15,-3){\line(0,-1){6}}
	\put(29,-3){\line(0,-1){6}}
	\put(18,-9){\line(0,-1){7.5}}
	\put(26,-9){\line(0,-1){7.5}}
	\put(20,-8){\footnotesize $\sigma$}
	\put(50,0){$+$}}
\hspace{23mm} \dfrac{1}{4!} \sum\limits_{\sigma \in S_4}
\mbox{\put(0,0){$\Big( \, \dfrac{b_2^2}{2}$}
	\put(45,3){\circle{40}}
	\put(38,5.5){\oval(11,9)}
	\put(52,5.5){\oval(11,9)}
	\put(36.5,1){\line(0,-1){6}}
	\put(40,1){\line(0,-1){6}}
	\put(50,1){\line(0,-1){6}}
	\put(53.5,1){\line(0,-1){6}}
	\put(35,-5){\line(1,0){20}}
	\put(35,-5){\line(0,-1){6}}
	\put(55,-5){\line(0,-1){6}}
	\put(35,-11){\line(1,0){20}}
	\put(43,-10){\footnotesize $\sigma$}
	\put(39.5,-11){\line(0,-1){5.5}}
	\put(43,-11){\line(0,-1){6}}
	\put(46,-11){\line(0,-1){6}}
	\put(49.5,-11){\line(0,-1){5.5}}
	\put(70,0){$+ \ \dfrac{b_4}{4!}$}
	\put(115,3){\circle{40}}
	\put(115,8){\oval(20,14)}
	\put(110,1){\line(0,-1){6}}
	\put(113,1){\line(0,-1){6}}
	\put(116,1){\line(0,-1){6}}
	\put(119.5,1){\line(0,-1){6}}
	\put(105,-5){\line(1,0){20}}
	\put(105,-5){\line(0,-1){6}}
	\put(125,-5){\line(0,-1){6}}
	\put(105,-11){\line(1,0){20}}
	\put(113,-10){\footnotesize $\sigma$}
	\put(110,-11){\line(0,-1){5.5}}
	\put(113,-11){\line(0,-1){6}}
	\put(116,-11){\line(0,-1){6}}
	\put(119.5,-11){\line(0,-1){5.5}}
	\put(138,0){$\Big) \ + \ \ldots$}}
\hspace{6.2cm}
\end{align}

Let us explicitly write down the second term:

\begin{align}
\Phi_{L}^{\text{adj}} \circ \rho \left( w_2 \right) = 
\dfrac{1}{2}\Big( \mbox{
	\put(25,3){\circle{40}}
	\put(25,4){\oval(20,14)}
	\put(21,-3){\line(0,-1){13.5}}
	\put(29,-3){\line(0,-1){13.5}}
	\put(50,0){$+$}
	\put(85,3){\circle{40}}
	\put(85,4){\oval(20,14)}
	\put(81,-3){\line(1,-2){7}}
	\put(89,-3){\line(-1,-2){7}}
	\put(110,0){\Big) \ =  }
	\put(156,3){\circle{40}}
	\put(156,4){\oval(20,14)}
	\put(152,-3){\line(0,-1){13.5}}
	\put(160,-3){\line(0,-1){13.5}}
	\put(180,0){$= \ (2t)^2 \cdot \dim L.$}}
\hspace{90mm}
\end{align}
\vspace{1mm}

Using IHX relation or just scrutinizing one of the diagrams it is easy to see that both terms of the second order are same. Terms of higher orders have much more complicated structure. Explicit calculation of the $\rho\left(w_4\right)$ can be found in Section 5.7 in the textbook \cite{chmutov2012introduction}. We list answers in the Table \ref{tab:symm_cas} for all terms up to 8 order, which we calculated with the help of the computer.
\renewcommand{\arraystretch}{2.4}
\begin{table}[h!]
\centering
\begin{tabular}{|c|c|}
	\hline
	Wheel diagram & $\phantom{\Big|}$ Value of $\Phi_{L}^{\text{adj}} \circ \rho$ $\phantom{\Big|}$  \\[4pt]
	\hline
	$w_2$ & $\phantom{\Big|}$ $(2t)^2 \cdot \dim L$ $\phantom{\Big|}$ \\[4pt]
	\hline
	$w_4$ & $\phantom{\Big|}$ $ \left( \dfrac{20}{3}\,t^4 - 3t\omega \right) \cdot \dim L $ $\phantom{\Big|}$ \\[4pt]
	\hline
	$w_2^2$ & $\phantom{\Big|}$ $ \dfrac{40}{3}t^4 \cdot \dim L $ $\phantom{\Big|}$ \\[4pt]
	\hline
	$w_6$ & $\phantom{\Big|}$ $ \left( \dfrac{28}{3}\,t^6 - \dfrac{42}{5}\,t^3\omega + \dfrac{15}{8}\,t\sigma\omega \right) \cdot \dim L $ $\phantom{\Big|}$ \\[4pt]
	\hline
	$w_2w_4$ & $\phantom{\Big|}$ $\left( \dfrac{56}{3}\,t^6 - \dfrac{42}{5}\,t^3\omega  \right) \cdot \dim L  $ $\phantom{\Big|}$ \\[4pt]
	\hline
	$w_2^3$ & $\phantom{\Big|}$ $ \dfrac{112}{3}\,t^6  \cdot \dim L $ $\phantom{\Big|}$ \\[4pt]
	\hline
	$w_8$ & $\phantom{\Big|}$ $ \left( 12t^8 - \dfrac{2867}{180}\,t^5\omega + \dfrac{51739}{10080}\,t^3\sigma\omega + \dfrac{59}{32}\,t^2\omega^2 -  \dfrac{161}{160}\,t\sigma^2\omega \right) \cdot \dim L $ $\phantom{\Big|}$ \\[4pt]
	\hline
	$w_2w_6$ & $\phantom{\Big|}$ $ \left( 24t^8 - \dfrac{109}{5}\,t^5\omega + \dfrac{1063}{280}\,t^3\sigma\omega + \dfrac{9}{8}\,t^2\omega^2  \right) \cdot \dim L $ $\phantom{\Big|}$ \\[4pt]
	\hline
	$w_4^2$ & $\phantom{\Big|}$ $ \left( 24t^8 - \dfrac{2474}{105}\,t^5\omega + \dfrac{239}{140}\,t^3\sigma\omega + \dfrac{141}{28}\,t^2\omega^2 -  \dfrac{9}{20}\,t\sigma^2\omega \right) \cdot \dim L $ $\phantom{\Big|}$ \\[4pt]
	\hline
	$w_2^2w_4$ & $\phantom{\Big|}$ $ \left( 48t^8 - \dfrac{364}{15}\,t^5\omega + \dfrac{6}{5}\,t^3\sigma\omega \right) \cdot \dim L $ $\phantom{\Big|}$ \\[4pt]
	\hline
	$w_2^4$ & $\phantom{\Big|}$ $ \left( 96t^8 - \dfrac{16}{5}\,t^5\omega \right) \cdot \dim L $ $\phantom{\Big|}$ \\[4pt]
	\hline
\end{tabular}
\caption{Universal expanding coefficients for the unknot}
\label{tab:symm_cas}
\end{table}

\subsubsection{Symmetric Casimir invariants}
\label{sec:symm-Cas}

The higher Casimir invariants are defined by formula \eqref{eq:higher_Cas_def} and by diagrams \eqref{eq:higher_Cas_as_diag}. One can note that $\Phi_{L}^{\text{adj}} \circ \rho(w_{2n})$ is an eigenvalue of the completely symmetrized $2n$-th Casimir invariant taken in the adjoint representation. Indeed, for usual Casimir invariant \eqref{eq:higher_Cas_as_diag} of degree $n$ we have 
\vspace{2mm}
\begin{align}
\text{tr}_{\text{adj}} ( C_n ) = (-1)^n \cdot 
\put(15,10){\circle{20}}
\put(7,4){\line(-1,-2){7}}
\put(9,2){\line(-1,-3){4}}
\put(11,1){\line(-1,-5){2.2}}
\put(11,-7){$...$}
\put(13,-4){\tiny{n}}
\put(19,1){\line(1,-5){2.2}}
\put(21,2){\line(1,-3){4}}
\put(23,4){\line(1,-2){7}}
\put(15,-13.3){\oval(45,7)}
\put(33,0){$= (-1)^n \cdot 2t \cdot \chi_{L}(\hat x_{n-1}) \cdot$}
\put(145,3){\circle{15}}
\put(160,0){$= (-1)^n \, 2t \, \chi_{L}(\hat x_{n-1}) \cdot \dim L$.}
\hspace{9cm}
\end{align}
In a similar way for the completely symmetrized Casimir invariants
\begin{align}\label{eq:symmCas-def}
C^s_n = \dfrac{1}{n!} \sum_{\sigma \in S_n} \Tr_{ad}(X^{\mu_1}X^{\mu_2}\cdots X^{\mu_p}) X_{\mu_1}X_{\mu_2}\cdots X_{\mu_p}
\end{align}
one can obtain
\begin{align}
\text{tr}_{\text{adj}} ( C^s_{2n} ) = \sum_{\sigma \in S_{2n}} \cdot 
\mbox{\put(25,20){\circle{20}}
\put(17,14){\line(-1,-2){5.5}}
\put(19,12){\line(-1,-3){3}}
\put(21,11){\line(-1,-5){1.6}}
\put(29,11){\line(1,-5){1.6}}
\put(31,12){\line(1,-3){3}}
\put(33,14){\line(1,-2){5.5}}
\put(21,5.5){...}
\put(21,-7){$...$}
\put(22,-1.5){\tiny{$\sigma$}}
\put(6,3){\line(1,0){38}}
\put(6,-3){\line(1,0){38}}
\put(6,3){\line(0,-1){6}}
\put(44,3){\line(0,-1){6}}
\put(9,-3){\line(-1,-2){3.9}}
\put(14,-3){\line(-1,-3){2.6}}
\put(18,-3){\line(-1,-5){1.6}}
\put(32,-3){\line(1,-5){1.6}}
\put(36,-3){\line(1,-3){2.6}}
\put(41.5,-3){\line(1,-2){3.9}}
\put(25,-14.3){\oval(50,7)}
\put(52,0){$= \ \Phi_{L}^{\text{adj}} \circ \rho \left( w_{2n} \right)$}
}
\hspace{4.5cm}
\end{align}

Such symmetric Casimir invariants were considered by Gelfand \cite{gel1950center}. They, just like the usual Casimir invariants, form a basis at the center of the universal enveloping space of $L$. In the paper \cite{mkrtchyan2012casimir} the fourth symmetric Casimir invariant was calculated.

We have not yet been able to find explicit expressions (or generating function) for all symmetric Casimir invariants $C^s_n$. However, it should be noted that the quantum dimension is a generating function for Casimir invariants of a certain type (defined by formula \eqref{eq:ZunknB}), which are generally no worse than Casimir invariants $C_n$ \eqref{eq:higher_Cas_def} or their symmetric versions $C^s_n$ \eqref{eq:symmCas-def}.

\subsection{Quantum knot polynomial}
\label{sec:quantum-knot-polynomial}

Quantum knot polynomials (i.e. the Jones, HOMFLY, Kauffman) are generalized by the Kontsevich knot invariant $ \text{I}: \mathcal{K} \rightarrow \mathcal{A} $ taking values in the space of Jacobi diagrams. Thus, it is very natural for it to admit Vogel universality. However, one should specify a mapping from \( \mathcal{A}\) to \( \Lambda\) and there are multiple ways of doing that.

The obvious one is to restrict the Kontsevich invariant to the adjoint representation, utilizing the map $ \pi_{\text{adj}}: \mathcal{A} \rightarrow \Gamma $ from the Jacobi diagrams to the space of 3-graphs, given by \eqref{eq:Wilson-to-adjoint}. We refer to this as the "adjoint Kontsevich invariant" \( \text{I}_{\text{adj}} := \pi_{\text{adj}} \circ \text{I}: \mathcal{K} \rightarrow \Gamma \).

Remember that $ \Gamma $ is "almost isomorphic" to $ \Lambda $. More accurately, $ \Gamma = \Gamma_0 \oplus \Gamma_{n \geq 1} = \langle \bigcirc\rangle \oplus \Lambda \cdot \thg$.
Remember the mapping $\Xi $ from $ \Gamma_{n \geq 1} $ to the Vogel's $ \Lambda $, given in \eqref{eq:3graphs-Lambda}. It is convenient to extend its action to the whole \( \Gamma \) by setting \( \Xi \left( \bigcirc \right)=0\). Then the composition $ \Xi \circ \text{I}_{\text{adj}}$ is a knot invariant taking values in $ \Lambda $, which we call universal knot polynomial.

This construction was considered by Patureau-Mirand in \cite{patureau2006quantum}, but the paper is more focused on $ D(2, 1, \alpha) $ knot polynomials. Other papers on this issue are \cite{mironov2016universal, bishler2025torus, mironov2025torus}, which aim to calculate the universal knot polynomials for torus knots based on the  Rosso-Jones formula. Note the difference between \( \Phi_L (\text{I}_{\text{adj}}) \) and \( \chi_L (\Xi \circ \text{I}_{\text{adj}})\). The former is the quantum knot polynomial in the adjoint representation \( \Phi^{\text{adj}}_{{L}} \circ \text{I}\), which is the most straightforward target for universalization and was considered in \cite{mironov2016universal}. The latter has the benefit of yielding nontrivial knot invariant for the \(D(2,1,\lambda)\) Lie superalgebra. They are related by \( \Phi_L (\text{I}_{\text{adj}}) =\text{dim}_L +  2 t \text{dim}_L \,\chi_L (\Xi \circ \text{I}_{\text{adj}}) \).

There are at least two known ways to define the Kontsevich invariant: via integration and the combinatorial. Ultimately, both approaches are perturbative, but the latter  seems to be more suitable for the the Vogel universality. 

In this construction, the Kontsevich invariant is made of only two building blocks: the Drinfeld associator and the exponent of a chord.
The procedure is to represent a knot as a closed braid. Then the  $\exp(\mbox{\tikz[scale=0.25, baseline=2pt]{\draw[thick, arrows={-{Stealth[length=1mm]}}] (0,0)--(0,1.5);\draw[thick, arrows={-{Stealth[length=1mm]}}] (1,0)--(1,1.5);\draw (0,0.59)--(1,0.59); }}/2)$ is assigned to each crossing, and the $ \Phi_{KZ} $ is used whenever one strand approaches the other one. A detailed review of this construction can be found in \cite{chmutov2012introduction}.

The $\exp(\mbox{\tikz[scale=0.25, baseline=2pt]{\draw[thick, arrows={-{Stealth[length=1mm]}}] (0,0)--(0,1.5);\draw[thick, arrows={-{Stealth[length=1mm]}}] (1,0)--(1,1.5);\draw (0,0.59)--(1,0.59); }}/2)$ in the adjoint representation is simply the \( \exp(-{\Psi}/2)\). The Drinfeld associator is a power series in $A$ and $B$, where $A$ and $B$ are the following two diagrams:

$$
\mbox{\tikz[scale=0.6]{
\coordinate (a) at (0,0);
\coordinate (b) at (5,0);
\draw[ultra thick, arrows={-{Stealth[length=3mm]}}] (a) -- +(0,2);
\draw[ultra thick, arrows={-{Stealth[length=3mm]}}] (a)++(1,0) -- +(0,2);
\draw[ultra thick, arrows={-{Stealth[length=3mm]}}] (a)++(2,0) -- +(0,2);
\draw[fill] (a)++(0,1) circle(0.07) -- +(1,0) circle(0.07);
\node at ($(a)+(-1,1)$) {\large \( A = \)};
\draw[ultra thick, arrows={-{Stealth[length=3mm]}}] (b) -- +(0,2);
\draw[ultra thick, arrows={-{Stealth[length=3mm]}}] (b)++(1,0) -- +(0,2);
\draw[ultra thick, arrows={-{Stealth[length=3mm]}}] (b)++(2,0) -- +(0,2);
\draw[fill] (b)++(1,1) circle(0.07) -- +(1,0) circle(0.07);
\node at ($(b)+(-1,1)$) {\large \( B = \)};
}
}
$$
The explicit formula for the Drinfeld associator was found by Le and Murakami \cite{le1996kontsevich}:

\begin{equation}\label{eq:LeMurakami}
\Phi_{KZ} = 1 + \sum_{m=2}^\infty \sum_{\substack{0<\mathbf{p}, \, 0<\mathbf{q},\\ |\mathbf{p}| + |\mathbf{q}| = m}} \eta(\mathbf{p,q}) \sum_{\substack{0\leq \mathbf{r} \leq \mathbf{p} \\ 0 \leq \mathbf{s} \leq \mathbf{q}}} (-1)^{|\mathbf{r}|+|\mathbf{j}|} \binom{\mathbf{p}}{\mathbf{r}} \binom{\mathbf{q}}{\mathbf{s}} B^{|\mathbf{s}|} (A,B)^{(\mathbf{i,j})}A^{|\mathbf{r}|}
\end{equation}
where $(A,B)^{(\mathbf{i,j})} = A^{i_1} B^{j_1} \cdots A^{i_l} B^{j_l}$, the \( \eta(\mathbf{p,q}) \) is related to the multiple zeta function.

The Drinfeld associator can be mapped to the adjoint representation by \( \pi_{\text{adj}}\), yielding an infinite sum over 6--legged diagrams. However, at the current stage we are not capable of operating with diagrams of high orders. 

On the Lie algebra side, the Drinfeld associator plays the role of the Racah matrices relating the $R$-matrices acting on different tensor factors (different strands of the braid) to each other by conjugation. In this direction there has been a development in \cite{mironov2016universalA}.

Interestingly, there exist other ways to map Jacobi diagrams \( \mathcal{A} \) to  \( \Lambda \). One of then is to consider the space of the open Jacobi diagrams \( \mathcal{B}\), which has  an advantage of the grading by the number of legs. This allows for a projection onto the subspace of the $2$-legged open Jacobi diagrams \(\mathcal{B}_2\), which is almost the same as \( \Lambda \). Precisely,  \( \mathcal{B}_2 = \langle-\rangle \oplus \Lambda \cdot \bbl \cong \Gamma \), which could be mapped to \( \Lambda \). This construction yields a different knot invariant with values in \( \Lambda\). 
This invariant might be of interest for some further research, because it circumvents the restriction to the adjoint sector in \( \text{I}_{\text{adj}}\), and also it incorporates the \( D_{17} \) diagram \ref{eq:kernel-examples}, undetected by the weight systems of all Lie (super)algebras. 

It is an interesting question whether an alternative non-perturbative definition of the universal knot polynomial could be given. Possibly by some skein relation or by an analogue of the Reshetikhin-Turaev approach via the R-matrix. However, it has yet to be developed. In the following subsection we briefly discuss how to perturbatively write the Kontsevich integral in a universal form for torus knots.
\bigskip

\subsubsection{Torus knots}
\label{sec:torus-knots}

Let us calculate the Kontsevich integral of the simplest torus knots up to 3rd order. There exists a formula similar to the formula for quantum dimension \eqref{eq:ZunknB}. In papers \cite{le1997parallel, bar2003two} there was derived the following expression for the Kontsevich integral for all torus knots $T[m,n]$:

\begin{align} \label{eq:KItorus}
	\text{I}(T[m,n]) = \psi^m \Big( \text{I}(\bigcirc) \texttt{\#} \exp\left( \dfrac{n}{2m} \bigcirc \put(-11.5,2.5){\line(1,0){8.5}} {\linethickness{0.3mm} \put(-7.2,2.5){\circle{9.5}}}\right) \Big),
\end{align}
where $\bigcirc$ is the zero-framed unknot{\iffalse\footnote{It means that the Kontsevich integral of the unknot include terms with isolated chords, while the resulting answer not.}\fi}, $\texttt{\#}$ is a connected sum of chord diagrams,\hspace{5.4pt} $\put(-5.3,2.5){\line(1,0){8.5}} {\linethickness{0.4mm} \put(-1.0,2.5){\circle{9.5}}}$ \hspace{1pt} is the chord diagram with one chord and $\psi$ is a special operation on Jacobi diagrams. For a chord diagram $D$ operation $\psi^m$ is defined as the sum of  all possible ways of lifting the ends of the chords to the \(m\)-sheeted connected covering of the Wilson loop of $D$:
\begin{align}
\psi^m 
\begin{picture}(150,20)(0,0)
\Big( 
{\linethickness{0.4mm}
	\put(20,3){\circle{30}}}
\put(5,3){\line(1,0){30}}
\put(5,3){\circle*{2.5}}
\put(35,3){\circle*{2.5}}
\put(40,0){$\Big) = \sum\limits_i \sum\limits_j $}
{\linethickness{0.4mm}
	\put(105,3){\circle{40}}
	\put(105,3){\circle{34}}
	\put(105,3){\circle{28}}
	\put(105,3){\circle{16}}}
\put(96.5,3){\line(1,0){17}}
\put(105,20){\color{white}\circle*{8}}
\put(111,21){\color{white}\circle*{8}}
\put(109,18){\color{white}\circle*{8}}
\put(105,16){\color{white}\circle*{8}}
\put(104,14){\color{white}\circle*{9}}
\put(101,13){\color{white}\circle*{7}}
\put(108,11){\color{white}\circle*{6}}
\put(109,10){\color{white}\circle*{6}}
\put(100,9){\color{white}\circle*{5}}
{\linethickness{0.4mm}
	\qbezier(101.9,22.8)(105,23)(114,17.4)
	\qbezier(101,19.6)(105,19.5)(112,15.2)
	\qbezier(98.1,15.3)(104,18)(111.4,7.8)
	\qbezier(98.1,7)(100.1,11.5)(103.2,14.5)
	\put(104.3,15.7){\line(1,1){1.6}}
	\put(107.4,18.4){\line(1,1){1.3}}
	\qbezier(110.1,20.7)(113,21.5)(115,20.4)
	\put(97,3){\circle*{1.5}}
	\put(113,3){\circle*{1.5}}
	\put(85,3){\circle*{1.5}}
	\put(88,3){\circle*{1.5}}
	\put(91,3){\circle*{1.5}}
	\multiput(85.5,3)(1.8,0){6}{\line(1,0){1}}
	\multiput(114.5,3)(1.8,0){6}{\line(1,0){1}}
	\put(119,3){\circle*{1.5}}
	\put(122,3){\circle*{1.5}}
	\put(125,3){\circle*{1.5}}}
\put(99,4){\tiny{$i$}}
\put(108,4){\tiny{$j$}}
\end{picture}
\end{align}
\vspace{1mm}

The simplest series of torus knots are 2-strand knots $T[2,n]$. The $\psi^2$ operation for chord diagrams up to and including order 3 can be applied relatively easily
\begin{align}
\psi^2 
\begin{picture}(300,20)(0,0)
\Big( 
\mbox{{\linethickness{0.4mm}
		\put(12,3){\circle{26}} }
	\put(-1,3){\line(1,0){26}}
	\put(-1,3){\circle*{2.5}}
	\put(25,3){\circle*{2.5}}
	\put(27,0){\Big) = 4} 
	{\linethickness{0.4mm}
		\put(71,3){\circle{26}} }
	\put(58,3){\line(1,0){26}}
	\put(84,3){\circle*{2.5}}
	\put(58,3){\circle*{2.5}}}
	\end{picture}
\end{align}
\vspace*{-5mm}
\begin{align}
\psi^2 
\begin{picture}(300,20)(0,0)
\Big( 
{\linethickness{0.4mm}
	\put(15,3){\circle{26}} }
\put(3,8){\line(1,0){24}}
\put(3,-2){\line(1,0){24}}
\put(3,8){\circle*{2.5}}
\put(27,8){\circle*{2.5}}
\put(3,-2){\circle*{2.5}}
\put(27,-2){\circle*{2.5}}
\put(30,0){\Big) = 12}
{\linethickness{0.4mm} 
	\put(76,3){\circle{26}} }
\put(64,8){\line(1,0){24}}
\put(64,-2){\line(1,0){24}}
\put(64,8){\circle*{2.5}}
\put(88,8){\circle*{2.5}}
\put(64,-2){\circle*{2.5}}
\put(88,-2){\circle*{2.5}}
\put(95,0){+ 4} 
{\linethickness{0.4mm}
	\put(127,3){\circle{26}} }
\put(118,12){\line(1,-1){18}}
\put(118,-6){\line(1,1){18}}
\put(118,12){\circle*{2.5}}
\put(136,12){\circle*{2.5}}
\put(118,-6){\circle*{2.5}}
\put(136,-6){\circle*{2.5}}
\end{picture}
\end{align}
\begin{align}
\psi^2 
\begin{picture}(300,20)(0,0)
\Big( 
{\linethickness{0.4mm}
	\put(15,3){\circle{26}} }
\put(6,12){\line(1,-1){18}}
\put(6,-6){\line(1,1){18}}
\put(6,12){\circle*{2.5}}
\put(24,12){\circle*{2.5}}
\put(6,-6){\circle*{2.5}}
\put(24,-6){\circle*{2.5}}
\put(30,0){\Big) =  8} 
{\linethickness{0.4mm}
	\put(76,3){\circle{26}} }
\put(64,8){\line(1,0){24}}
\put(64,-2){\line(1,0){24}}
\put(64,8){\circle*{2.5}}
\put(88,8){\circle*{2.5}}
\put(64,-2){\circle*{2.5}}
\put(88,-2){\circle*{2.5}}
\put(95,0){+ 8} 
{\linethickness{0.4mm}
	\put(127,3){\circle{26}} }
\put(118,12){\line(1,-1){18}}
\put(118,-6){\line(1,1){18}}
\put(118,12){\circle*{2.5}}
\put(136,12){\circle*{2.5}}
\put(118,-6){\circle*{2.5}}
\put(136,-6){\circle*{2.5}}
\end{picture}
\end{align}
\vspace*{-5mm}
\begin{align}
	\psi^2 
\mbox{\begin{picture}(300,20)(0,0)
\Big( 
{\linethickness{0.4mm}
	\put(15,3){\circle{26}} }
\put(4,10){\line(1,0){22}}
\put(2,3){\line(1,0){26}}
\put(4,-4){\line(1,0){22}}
\put(4,10){\circle*{2.5}}
\put(26,10){\circle*{2.5}}
\put(2,3){\circle*{2.5}}
\put(28,3){\circle*{2.5}}
\put(4,-4){\circle*{2.5}}
\put(26,-4){\circle*{2.5}}
\put(30,0){\Big) = 32}
{\linethickness{0.4mm} 
	\put(76,3){\circle{26}} }
\put(65,10){\line(1,0){22}}
\put(63,3){\line(1,0){26}}
\put(65,-4){\line(1,0){22}}
\put(65,10){\circle*{2.5}}
\put(87,10){\circle*{2.5}}
\put(63,3){\circle*{2.5}}
\put(89,3){\circle*{2.5}}
\put(65,-4){\circle*{2.5}}
\put(87,-4){\circle*{2.5}}
\put(95,0){+ 16} 
{\linethickness{0.4mm}
	\put(132,3){\circle{26}} }
\put(121,10){\line(1,0){22}}
\put(119,3){\line(7,-2){24}}
\put(121,-4){\line(7,2){24}}
\put(121,10){\circle*{2.5}}
\put(143,10){\circle*{2.5}}
\put(119,3){\circle*{2.5}}
\put(145,3){\circle*{2.5}}
\put(121,-4){\circle*{2.5}}
\put(143,-4){\circle*{2.5}}
\put(151,0){+ 16} 
{\linethickness{0.4mm}
	\put(188,3){\circle{26}} }
\put(183,15){\line(0,-1){24}}
\put(193,15){\line(0,-1){24}}
\put(175,3){\line(1,0){26}}
\put(183,15){\circle*{2.5}}
\put(193,15){\circle*{2.5}}
\put(175,3){\circle*{2.5}}
\put(201,3){\circle*{2.5}}
\put(183,-9){\circle*{2.5}}
\put(193,-9){\circle*{2.5}}
\end{picture}}
\end{align}

\vspace*{-7mm}
\begin{align}
	       \psi^2 
\mbox{\begin{picture}(300,25)(0,0)
\Big( 
{\linethickness{0.4mm}
	\put(15,3){\circle{26}} }
\put(4,10){\line(1,0){22}}
\put(2,3){\line(7,-2){24}}
\put(4,-4){\line(7,2){24}}
\put(4,10){\circle*{2.5}}
\put(26,10){\circle*{2.5}}
\put(2,3){\circle*{2.5}}
\put(28,3){\circle*{2.5}}
\put(4,-4){\circle*{2.5}}
\put(26,-4){\circle*{2.5}}
\put(30,0){\Big) = 20}
{\linethickness{0.4mm} 
	\put(76,3){\circle{26}} }
\put(65,10){\line(1,0){22}}
\put(63,3){\line(1,0){26}}
\put(65,-4){\line(1,0){22}}
\put(65,10){\circle*{2.5}}
\put(87,10){\circle*{2.5}}
\put(63,3){\circle*{2.5}}
\put(89,3){\circle*{2.5}}
\put(65,-4){\circle*{2.5}}
\put(87,-4){\circle*{2.5}}
\put(95,0){+ 28} 
{\linethickness{0.4mm}
	\put(132,3){\circle{26}} }
\put(121,10){\line(1,0){22}}
\put(119,3){\line(7,-2){24}}
\put(121,-4){\line(7,2){24}}
\put(121,10){\circle*{2.5}}
\put(143,10){\circle*{2.5}}
\put(119,3){\circle*{2.5}}
\put(145,3){\circle*{2.5}}
\put(121,-4){\circle*{2.5}}
\put(143,-4){\circle*{2.5}}
\put(151,0){+ 12} 
{\linethickness{0.4mm}
	\put(188,3){\circle{26}} }
\put(183,15){\line(0,-1){24}}
\put(193,15){\line(0,-1){24}}
\put(175,3){\line(1,0){26}}
\put(183,15){\circle*{2.5}}
\put(193,15){\circle*{2.5}}
\put(175,3){\circle*{2.5}}
\put(201,3){\circle*{2.5}}
\put(183,-9){\circle*{2.5}}
\put(193,-9){\circle*{2.5}}
\put(205,0){+ 4} 
{\linethickness{0.4mm}
	\put(238,3){\circle{26}} }
\put(233,15){\line(2,-5){9.5}}
\put(243,15){\line(-2,-5){9.5}}
\put(225,3){\line(1,0){26}}
\put(233,15){\circle*{2.5}}
\put(243,15){\circle*{2.5}}
\put(225,3){\circle*{2.5}}
\put(251,3){\circle*{2.5}}
\put(233.4,-9){\circle*{2.5}}
\put(242.6,-9){\circle*{2.5}}
\end{picture}}
\end{align}

\vspace*{3mm}
The explicit expression for $\text{I}(\bigcirc)$ is given in Section \ref{sec:qdim} above. Then we can calculate the explicit expression of formula \eqref{eq:KItorus} up to 3 order for $T[2,n]$:

\begin{equation}
	\mbox{\tikz[scale=0.46, baseline = -3pt]{
			\coordinate (a) at (-1,0);
			\coordinate (b) at (3,0);
			\coordinate (c) at (7,0);
			\coordinate (d) at (12.4,0);
			\coordinate (e) at (16.5,0);
			\coordinate (f) at (22.4,0);
			\coordinate (g) at (27.5,0);
			
			\node at ($(a)+(-3.2,0)$) {\( \text{I} (T[2,n]) = \)};
			\node at ($(b)+(-1.8, 0)$) {\( + \, n \)};
			\node at ($(c)+(-1.9, 0)$) {\( + \, \frac{ n^2 \hspace{-1mm}}{2} \,\)};
			\node at ($(d)+(-2.5, 0)$) {\( - \, \frac{3n^2-4}{48} \)};
			\node at ($(e)+(-1.9, 0)$) {\( + \, \frac{ n^3 \hspace{-1mm}}{6} \,\)};
			
			\node at ($(f)+(-2.75, 0)$) {\( - \, \frac{3n^3-4n}{48} \)};
			\node at ($(g)+(-2.45,0)$) {\( + \; \frac{n^3-n}{96} \)};

			\draw[ultra thick] (a) circle(1);
			\draw[ultra thick] (b) circle(1);
			\draw[ultra thick] (c) circle(1);
			\draw[ultra thick] (d) circle(1);
			\draw[ultra thick] (e) circle(1);
			\draw[ultra thick] (f) circle(1);
			\draw[ultra thick] (g) circle(1);
			
			\draw (b)+(180:1) -- +(0:1);
			
			\draw (c)+(160:1) -- +(20:1);
			\draw (c)+(200:1) -- +(-20:1);
			
			\draw (d) circle (0.3);
			\draw (d)+(0:0.3) -- +(0:1);
			\draw (d)+(180:0.3) -- +(180:1);
			
			\draw (e)+(180:1) -- +(0:1);
			\draw (e)+(150:1) -- +(30:1);
			\draw (e)+(210:1) -- +(-30:1);

			\draw (f)+(0, -0.25) circle(0.2);
			\draw (f)+(-0.2, -0.25) -- +(194.5:1);
			\draw (f)+(0.2, -0.25) -- +(-14.5:1);
			\draw (f)+(155:1) -- +(25:1);
			
			\draw (g) +(-0.35, 0) circle (0.2);
			\draw (g) +(0.35, 0) circle (0.2);
			\draw (g)+(180:1) -- +(180:0.55);
			\draw (g)+(180:0.15) -- +(0:0.15);
			\draw (g)+(0:1) -- +(0:0.55);
	}}
\end{equation}

Let us stress that $\text{I} (T[2,n])$ is a framed version of the Kontsevich integral. To obtain the unframed version one needs to impose 1T relation, which means that all Jacobi diagrams containing an isolated chord are set to zero (see the textbook \cite{chmutov2012introduction} for details of deframing procedure):

\begin{equation}
	\mbox{\tikz[scale=0.46, baseline = -3pt]{
			\coordinate (a) at (-1,0);
			\coordinate (d) at (4.4,0);
			\coordinate (g) at (9.6,0);
			
			\node at ($(a)+(-3.8,0)$) {\( \text{I}^{\text{unfr}} (T[2,n]) = \)};
			\node at ($(d)+(-2.5, 0)$) {\( - \, \frac{3n^2-4}{48} \)};
			
			\node at ($(g)+(-2.5,0)$) {\( + \; \frac{n^3-n}{96} \)};

			\draw[ultra thick] (a) circle(1);
			\draw[ultra thick] (d) circle(1);
			\draw[ultra thick] (g) circle(1);

			\draw (d) circle (0.3);
			\draw (d)+(0:0.3) -- +(0:1);
			\draw (d)+(180:0.3) -- +(180:1);

			\draw (g) +(-0.35, 0) circle (0.2);
			\draw (g) +(0.35, 0) circle (0.2);
			\draw (g)+(180:1) -- +(180:0.55);
			\draw (g)+(180:0.15) -- +(0:0.15);
			\draw (g)+(0:1) -- +(0:0.55);
	}} \ .
\end{equation}

Therefore, the adjoint Kontsevich integral takes the following form:
\begin{align}
\Phi_{L}^{\text{adj}} \left( \text{I}^{\text{unfr}}(T[2,n]) \right) = \dim L \cdot \Big( 1 \ - \ \dfrac{3n^2-1}{12}\, t^2  \ + \ \dfrac{n^3+n}{12}\, t^3   \ + \ \ldots \Big).
\end{align}

\section{Universality beyond Lie algebras}
\label{sec:discussion}

In today's literature, Vogel's universality is almost always understood as the universal behavior of certain quantities in Lie algebras restricted to the adjoint sector of representations (for attempts to extend this to other representations, see \cite{suprun2026vogel, isaev2026vogel}). To the best of our knowledge, however, no explanation for the phenomenon of universality within the framework of Lie algebras has ever been proposed. Vogel himself originally started from the diagrammatic $\Lambda$-algebra, yet his approach has received virtually no further development over the past twenty years. This lack of progress can likely be attributed to the fact that the theory of Lie algebras and their representations is far more developed than the diagrammatic theory. Nevertheless, the diagrammatic approach offers several key \textit{advantages}. First, it at least partially explains why universal formulas for certain quantities can be written. Second, it allows such formulas to be derived directly in a universal form, eliminating the need to perform calculations for each Lie algebra individually and subsequently construct an analytic continuation. In this paper, we have attempted to highlight the benefits of adopting the diagrammatic approach. While this method remains underdeveloped and, as such, is not yet capable of performing complex calculations efficiently, we hope to have convincingly demonstrated its unique conceptual appeal.

In this section, we focus on two closely interrelated questions: What is universality? and What is its origin?

As noted above, the answer to the second question remains unknown. This is quite surprising, given that the history of Lie algebras spans over 150 years, and the phenomenon of universality itself was discovered more than three decades ago. Despite sustained attention from numerous distinguished specialists, including Pierre Deligne \cite{deligne1996serie, deligne2002exceptional}, no explanation in terms of Lie algebras has yet been proposed.

Answering the first question is, in fact, also nontrivial. If we say that the phenomenon of universality consists in the fact that certain quantities, such as the dimensions of representations, behave uniformly and can be parameterized by the three Vogel parameters $(\alpha,\beta,\gamma)$, then the following question immediately arises: In what class of functions are we considering this parameterization? Obviously, if we allow all possible functions, the phenomenon itself disappears, because describing a discrete set of quantities (like dimensions) is trivial using a tabular function of a single variable. The phenomenon becomes meaningful only if we greatly narrow the class of functions, for example, by requiring them to be rational functions. But then the question immediately follows: why should we restrict ourselves precisely to this class? Clearly, this question is closely related to the origin of the universality phenomenon.

At the same time, Vogel's original explanation of universality (via the characters of the $\Lambda$-algebra) must be regarded as, at the very least, incomplete. Consider the example of the dimension of the adjoint representation \eqref{eq:univ-dim}. The diagrammatic technique yields a formula expressed in variables $(t, \sigma, \omega)$, which have a natural interpretation in terms of diagrams. However, when we switch to the Vogel parameters $(\alpha,\beta,\gamma)$, the numerator and denominator of this formula factor into linear factors. From the perspective of diagrammatic techniques, it is no longer clear how to explain this factorization, since the variables $(\alpha,\beta,\gamma)$ lack a natural diagrammatic interpretation. This example reveals that for Lie algebra weight systems, the phenomenon of universality possesses a broader set of properties.

Based on the above considerations, we believe it is both important and useful within Vogel's diagrammatic approach to separate those cases that can be obtained at a purely diagrammatic level, without invoking the concept of a weight system.

\bigskip

$\textbf{1.}$ \textbf{Casimir operators.} Casimir operators acting in the adjoint representation are represented by 2-legged trivalent diagrams as we discussed in Section \ref{sec:Casimirs}. Formula \eqref{eq:higher_Cas_as_diag} 
\begin{equation}
	\mbox{\tikz[baseline = (current bounding box.center)]{
			\node at (-1.5, 0.75) {$(-1)^p C_p \ =$};
			\node at (0, 0.75) {$\hat{x}_{p-1} \ \cdot$};
			\draw (0.7, 0.75) -- +(0.5, 0);
			\draw (1.7, 0.75) circle (0.5);
			\draw (2.2, 0.75) -- (2.7, 0.75);
	}}
\end{equation}
shows that any operator $C_p$ in the adjoint representation is given by the multiplication of the corresponding element of $\Lambda$-algebra $\hat{x}_{p-1}$ and the 2-legged diagram $-\hspace{-1.5mm}\bigcirc\hspace{-1.5mm}-$. Therefore, one can represent all Casimir operators by 2-legged trivalent diagrams on which the Vogel $\Lambda$-algebra acts naturally.

Note that at the diagrammatic level, Casimir operators require trivalent diagrams of two types: two-legged and three-legged. If we consider Lie algebra weight system then we get \textbf{one} additional relation  $$\Phi_{L}\left(-\hspace{-1.5mm}\bigcirc\hspace{-1.5mm}-\right) = 2t \cdot \rm{Id}.$$

\

$\textbf{2.}$ \textbf{Projectors.} 
Consider a vector space spanned by 4-legged trivalent diagrams modulo AS+IHX relations. Define a multiplication by gluing from the right. Then let us consider the following diagrams:
\begin{equation}
	\mbox{\begin{picture}(-20,10)
			\put(-100,0){$P_a = \dfrac{1}{2} \Big($}
			\put(-60,-3){\line(1,0){20}}
			\put(-60,9){\line(1,0){20}}
			\put(-35,0){$-$}
			\put(-20,-3){\line(5,2){13}}
			\put(10,9){\line(-5,-2){13}}
			\put(-20,9){\line(5,-2){30}}
			\put(15,0){$\Big),$}
			\put(35,0){$P_a^2 = P_a$}
	\end{picture}}
\end{equation}

\begin{equation}
	\begin{picture}(-20,10)
		\put(-100,0){$P_s = \dfrac{1}{2} \Big( $}
		\put(-60,-3){\line(1,0){20}}
		\put(-60,9){\line(1,0){20}}
		\put(-35,0){$+$}
		\put(-20,-3){\line(5,2){13}}
		\put(10,9){\line(-5,-2){13}}
		\put(-20,9){\line(5,-2){30}}
		\put(15,0){$\Big),$}
		\put(35,0){$P_s^2 = P_s$}
	\end{picture}
\end{equation}

\begin{equation}
	\mbox{\begin{picture}(30,10)
			\put(-50,0){$P_{0} = $}
			\put(-20,3){\oval(14,14)[r]}
			\put(5,3){\oval(14,14)[l]}
			\put(8,0){$,$}
			\put(20,0){$P_{0}^2 = \bigcirc \cdot P_{0}$}
	\end{picture}}
\end{equation}

\begin{equation}
	\mbox{\begin{picture}(-20,10)
			\put(-94,0){$P_{\text{adj}} = $}
			\put(-60,9){\line(1,0){10}}
			\put(-60,-3){\line(1,0){10}}
			\put(-50,-3){\line(1,1){6}}
			\put(-50,9){\line(1,-1){6}}
			\put(-44,3){\line(1,0){10}}
			\put(-34,3){\line(1,1){6}}
			\put(-34,3){\line(1,-1){6}}
			\put(-28,9){\line(1,0){10}}
			\put(-28,-3){\line(1,0){10}}
			\put(-15,0){$,  \quad  P_{\text{adj}}^2 = 2\hat{t} \cdot P_{\text{adj}}$}
	\end{picture}}
\end{equation}

Note that $P_a$ and $P_s$ are projectors. The diagrams $P_0$ and $P_{\text{adj}}$ are not true projectors because of the factors $\bigcirc$ and $2\hat{t}$. To obtain genuine projectors, one should divide by these factors. To do so, it is necessary to extend the field over which the algebra of such diagrams is considered, i.e., adjoin the elements $\bigcirc$ and $\hat{t}$. While it is possible to extend the field by adjoining $\bigcirc$, \(\hat{t}\) being a zero divisor (see Section~\ref{sec:zero-divisor}) becomes an obstacle.

\

$\textbf{3.}$ \textbf{Kontsevich integral.} As we already discussed in Section \ref{sec:quantum-knot-polynomial} adjoint Kontsevich integral $\pi_{\text{adj}} \circ \text{I}(\mathcal{K})$ takes values in the algebra of 3-graphs, which is isomorphic to $\Lambda$-algebra (up to 0 degree element $\langle\bigcirc\rangle$). Thus, the adjoint Kontsevich integral can be expressed completely in terms of diagrams from the Vogel $\Lambda$-algebra.

\section*{Acknowledgements}
We are grateful for fruitful discussions to L.Bishler, A.Isaev, S.Krivonos, A.Mironov, R.Mkrtchayn, Al.Morozov, An.Morozov, and A.Provorov. We also are indebted for discussions to the organizers and participants of the workshop "Universal description of Lie algebras, Vogel theory, applications" in JINR, Dubna.

\noindent Funding for this publication was generously provided by the Priority 2030 Academic Leadership Initiative, contributing to the educational work of "Universities for a New Generation of Leaders", a project within the framework of the federal Youth and Children program.

\printbibliography

\end{document}